# The Morozov's principle applied to data assimilation problems

Laurent Bourgeois and Jérémi Dardé

June 7, 2022

**Abstract**

This paper is focused on the Morozov's principle applied to an abstract data assimilation framework, with particular attention to three simple examples: the data assimilation problem for the Laplace equation, the Cauchy problem for the Laplace equation and the data assimilation problem for the heat equation. Those ill-posed problems are regularized with the help of a mixed type formulation which is proved to be equivalent to a Tikhonov regularization applied to a well-chosen operator. The main issue is that such operator may not have a dense range, which makes it necessary to extend well-known results related to the Morozov's choice of the regularization parameter to that unusual situation. The solution which satisfies the Morozov's principle is computed with the help of the duality in optimization, possibly by forcing the solution to satisfy given a priori constraints. Some numerical results in two dimensions are proposed in the case of the data assimilation problem for the Laplace equation.

# 1 Introduction

In this article we consider inverse problems of the following type: find a solution $u$ to the Laplace equation $\Delta u = 0$ in a domain $\Omega$ of $\mathbb{R}^d$, $d \geq 1$, from values of $u$ measured on a smaller domain $\omega$. Such problem can be seen as a toy model for data assimilation problems. Data assimilation is a very active domain of applied mathematics in connection with oceanography, meteorology or life sciences. It can be addressed from a deterministic or a stochastic point of view, one of the pioneering contribution to the field being [2], a deterministic vision being exposed in [24], a stochastic one in [22]. The introduction of [14] also offers a nice overview of data assimilation problems. Despite uniqueness holds for our basic data assimilation problem for the Laplace equation, it is severely ill-posed, in the sense that existence is obviously not ensured if the measurements are corrupted by noise. In order to regularize such problem, we introduce a mixed variational formulation which is parametrized by a small parameter $\varepsilon > 0$. Such formulation has the advantage to be well-posed and to provide a solution which is close to the true solution $u$. In addition, the regularized solution is searched in the natural energy space in which $u$ is supposed to belong, typically $H^1(\Omega)$, which results from the mixed nature of our formulation. The "mixed" terminology comes from the fact that the weak formulation consists of a system of two equations satisfied by two unknowns, that is the regularized solution $u_\varepsilon$ and a Lagrange multiplier $\lambda_\varepsilon$, which is also searched in the space $H^1(\Omega)$. The idea of mixed formulation to solve linear ill-posed problems goes back to [3], which concerns the Cauchy problem for the Laplace equation and can be seen as a variant of the quasi-reversibility method introduced in [21]. The notion of mixed formulation was more recently recast in a general abstract framework in [8], in which it is proven that such mixed formulation is equivalent to the classical Tikhonov regularization for a well-chosen linear injective operator $A$. In particular, once the mixed formulation has been reformulated as a Tikhonov regularization, it is natural to apply the classical Morozov's discrepancy principle to compute a consistent regularization parameter $\varepsilon$ as a function of the amplitude of noise $\delta$. In [8], this technique was used in the case of a Cauchy-type problem for the Helmholtz equation set in a waveguide, that is the data was formed by the trace and the normal derivative of the acoustic field on a subpart of the boundary. It is important to note that for such a Cauchy-type problem, the underlying operator $A$ has a dense range.

When it comes to the data assimilation problem for the Laplace equation previously introduced, we will observe that applying the Morozov's discrepancy principle is not standard any more in the sense that the corresponding operator $A$ has not a dense range. Indeed, the standard results which justify the Morozov's principle for the Tikhonov regularization (see for example Theorem 2.17 in [20]) are limited to dense range operators. The objectives of the present paper are the following. The first objective is to generalize these standard results to the case when such denseness assumption fails. In particular, we will see that this generalized result requires the data to satisfy an additional condition which is not trivial to check. A second objective is to extend the duality method, introduced in the context of Morozov's discrepancy principle in [4] and revisited in [6] and [12], to the case of those operators $A$ which do not have a dense range. This duality method, the origin of which is the theory developed in [17], later adapted to controllability problems in [25] and in [26], consists in solving an unconstrained minimization problem involving a cost function which depends on the adjoint operator $A^*$ of $A$. The main interest of the duality method in the context of inverse problems consists of an idea introduced in [4, 6] and reused in the present paper: by applying the operator $A$ to the solution of such minimization problem, we exactly obtain the Tikhonov solution associated with the Morozov's choice for $\varepsilon$. In our paper, we also adapt another idea, borrowed from [18] in the context of control theory, then transposed in [12] in the context of inverse problems. It consists in introducing a modification of the cost function involving a compact projection operator in order to impose some *a priori* assumptions to the solution which may be useful in practice, while keeping the objective of satisfying the Morozov's principle. The third objective is to apply all the previous ideas to an abstract framework of data assimilation problems, including the toy problem presented at the beginning of the introduction.

Although we present some numerical results using a finite element method at the end of the paper, the choice of the discretization parameter $h$ is not discussed in the present paper, such $h$ being supposed to be sufficiently small so that we can apply the algorithms introduced at the continuous level. The Morozov's principle for the discretized Tikhonov regularization is however an interesting subject, addressed for example in [27]. In a long series of papers, Burman and his collaborators (see for example [10, 11, 9] for the Laplace, heat and wave equations, respectively) have proposed a weak formulation which is different from ours: there are no regularization parameter at the continuous level, but sophisticated stabilizers at the discrete level are introduced to obtain a well-posed problem. In those works, some estimates between the exact solution and the solution to the discrete regularized problem are obtained with respect to the mesh parameter $h$, which can be seen as a regularization parameter. Similar estimates are obtained for discretized mixed formulations to regularize the Cauchy problem for the Laplace equation in [5] and a data assimilation problem for a wave-type equation in [15]. Those estimates are interesting because they provide convergence rates, up to a multiplicative constant which cannot be estimated. The subject of our article is different and also challenging: we wish to choose a particular value of the regularization parameter.

Our paper is organized as follows. In section 2, we revisit the Morozov's principle for the Tikhonov regularized solutions in the case of an operator which has not a dense range. An abstract framework for data assimilation problems is introduced in section 3, to which the results obtained in section 2 are applied. In section 4, we analyze three particular cases of our general data assimilation framework: the data assimilation problem for the Laplace equation, the Cauchy problem for the Laplace equation, lastly the data assimilation problem for the heat equation. Section 5 is dedicated to numerical experiments related to the first particular case. Lastly, an appendix summarizes the theory of duality in optimization exposed in [17].

## 2 The Morozov's principle revisited

### 2.1 The Tikhonov regularization and the Morozov's principle

We first extend a well-known result related to the Tikhonov regularization for operators which have a dense range (see for example section 2.5 in [20] for the restricted case of a compact operator) to the case of operators which don't satisfy this property. Let $A : V \to H$ be a linear bounded operator from the Hilbert space $V$ to another Hilbert space $H$. We assume that $A$ is injective. In what follows, for a family $(u_\varepsilon)_\varepsilon$ of functions depending on the real parameter $\varepsilon > 0$, we will frequently use a slight abuse of notations: each time the family $(u_\varepsilon)_\varepsilon$ will be called a sequence, we mean that we may choose any sequence $(\varepsilon_n)_{n\in\mathbb{N}}$, such that $\varepsilon_n \to 0$ when $n \to +\infty$. Extracting a subsequence of $(u_\varepsilon)_\varepsilon$ means that we have extracted a subsequence of that particular sequence of $(\varepsilon_n)_{n\in\mathbb{N}}$.

**Theorem 2.1.** *For some $\delta > 0$, let us assume that the data $g^\delta \in H$ is such that*

$$\|g_\perp^\delta\|_H < \delta < \|g^\delta\|_H, \tag{1}$$

*where $g_\perp^\delta$ is the orthogonal projection of $g^\delta$ on $(\mathrm{Range}\, A)^\perp$.*
*For $\varepsilon > 0$, let us denote $u_\varepsilon^\delta \in V$ the regularized solution associated with data $g^\delta$ in the sense of Tikhonov, which is defined by the weak formulation*

$$(Au_\varepsilon^\delta, Av)_H + \varepsilon(u_\varepsilon^\delta, v)_V = (g^\delta, Av)_H, \quad \forall v \in V. \tag{2}$$

*There exists a unique $\varepsilon > 0$ such that*
$$\|Au_\varepsilon^\delta - g^\delta\|_H = \delta. \tag{3}$$

*Proof.* Let us introduce, for $\varepsilon > 0$, the function
$$E^\delta(\varepsilon) = \|Au_\varepsilon^\delta - g^\delta\|_H^2,$$
which is differentiable and satisfies
$$\frac{dE^\delta}{d\varepsilon}(\varepsilon) = 2(Au_\varepsilon^\delta - g^\delta, Av_\varepsilon^\delta)_H, \quad \forall \varepsilon > 0,$$
where $v_\varepsilon^\delta \in V$ is uniquely defined by
$$(Av_\varepsilon^\delta, Av)_H + \varepsilon(v_\varepsilon^\delta, v)_V = -(u_\varepsilon^\delta, v)_V, \quad \forall v \in V. \tag{4}$$
By choosing $v = v_\varepsilon^\delta$ in (2),
$$\frac{dE^\delta}{d\varepsilon}(\varepsilon) = -2\varepsilon(u_\varepsilon^\delta, v_\varepsilon^\delta)_V, \quad \forall \varepsilon > 0.$$
Then choosing $v = v_\varepsilon^\delta$ in (4), we find
$$\frac{dE^\delta}{d\varepsilon}(\varepsilon) = 2\varepsilon(\|Av_\varepsilon^\delta\|_H^2 + \varepsilon\|v_\varepsilon^\delta\|_V^2), \quad \forall \varepsilon > 0.$$

Obviously $dE^\delta/d\varepsilon \geq 0$, and more precisely $dE^\delta/d\varepsilon > 0$. Indeed, let us assume that $dE^\delta/d\varepsilon = 0$. It follows from the above identity that $v_\varepsilon^\delta = 0$, and from (4) that $u_\varepsilon^\delta = 0$. From (2) we hence infer that $(g^\delta, Av)_V = 0$ for all $v \in V$, that is $g^\delta \in (\text{Range } A)^\perp$, in other words $g^\delta = g_\perp^\delta$. But there is a contradiction with the fact that $\|g_\perp^\delta\|_H < \delta < \|g^\delta\|_H$. We conclude that the function $E^\delta$ is a continuous and non-decreasing function in the interval $(0, +\infty)$.

Now let us prove that
$$\lim_{\varepsilon \to 0} E^\delta(\varepsilon) = \|g_\perp^\delta\|_H^2 \quad \text{and} \quad \lim_{\varepsilon \to +\infty} E^\delta(\varepsilon) = \|g^\delta\|_H^2.$$
By choosing $v = u_\varepsilon^\delta$ in (2), we firstly obtain that
$$\|Au_\varepsilon^\delta\|_H \leq \|g^\delta\|_H, \tag{5}$$
and secondly that
$$\|u_\varepsilon^\delta\|_V \leq \frac{\|g^\delta\|_H}{\sqrt{\varepsilon}}. \tag{6}$$
This implies that $(u_\varepsilon^\delta)_\varepsilon$ converges to 0 in $V$ when $\varepsilon \to +\infty$, and then $E^\delta(\varepsilon)$ converges to $\|g^\delta\|_H^2$ when $\varepsilon$ tends to $+\infty$.

Since the function $E^\delta$ is non-negative and non-decreasing, it has a limit when $\varepsilon$ tends to 0. As by (5), $(Au_\varepsilon^\delta)_\varepsilon$ is bounded, we can extract from $(u_\varepsilon^\delta)_\varepsilon$ a subsequence, that we still denote $(u_\varepsilon^\delta)_\varepsilon$, such that $(Au_\varepsilon^\delta)_\varepsilon$ weakly converges to some $h \in H$. In fact, $h \in \overline{\text{Range } A}$ because the space $\overline{\text{Range } A}$ is weakly closed. We have, in view of (6),
$$|(Au_\varepsilon^\delta - g^\delta, Av)_H| = |-\varepsilon(u_\varepsilon^\delta, v)_V| \leq \sqrt{\varepsilon}\|g^\delta\|_H\|v\|_V \to 0, \quad \forall v \in V,$$
when $\varepsilon \to 0$. But we also have
$$(Au_\varepsilon^\delta - g^\delta, Av)_H \to (h - g^\delta, Av)_H, \quad \forall v \in V,$$
when $\varepsilon \to 0$. We conclude that $h - g^\delta \in (\text{Range } A)^\perp$. Since we have the decomposition $H = \overline{\text{Range } A} \bigoplus (\text{Range } A)^\perp$, we get that $h = g_{/\!/}^\delta$, where $g_{/\!/}^\delta$ is the orthogonal projection of data $g^\delta$ on $\overline{\text{Range } A}$, while $g^\delta - g_{/\!/}^\delta = g_\perp^\delta$. Hence $(Au_\varepsilon^\delta)_\varepsilon$ weakly converges to $g_{/\!/}^\delta$ in $H$, and subsequently, $(Au_\varepsilon^\delta - g^\delta)_\varepsilon$ weakly converges to $-g_\perp^\delta$ in $H$.

On the one hand, we have
$$\|g_\perp^\delta\|_H^2 \leq \liminf_{\varepsilon \to 0} \|Au_\varepsilon^\delta - g^\delta\|_H^2 = \liminf_{\varepsilon \to 0} E^\delta(\varepsilon) = \lim_{\varepsilon \to 0} E^\delta(\varepsilon).$$

On the other hand, we have

$$\begin{aligned}
E^\delta(\varepsilon) &= \|Au_\varepsilon^\delta - g^\delta\|_H^2 = (Au_\varepsilon^\delta - g^\delta, Au_\varepsilon^\delta)_H - (Au_\varepsilon^\delta - g^\delta, g^\delta)_H \\
&= -\varepsilon\|u_\varepsilon^\delta\|_V^2 - (Au_\varepsilon^\delta - g^\delta, g^\delta)_H \leq -(Au_\varepsilon^\delta - g^\delta, g^\delta)_H,
\end{aligned}$$

which implies that

$$\lim_{\varepsilon \to 0} E^\delta(\varepsilon) \leq (g_\perp^\delta, g^\delta)_H = \|g_\perp^\delta\|_H^2.$$

We conclude that $E^\delta(\varepsilon) \to \|g_\perp^\delta\|_H^2$ when $\varepsilon \to 0$. Finally since $E^\delta(0^+) = \|g_\perp^\delta\|_H^2 < \delta^2 < E^\delta(+\infty) = \|g^\delta\|_H^2$, since $E^\delta$ is a non-decreasing continuous function of $\varepsilon$, there exists a unique $\varepsilon > 0$ such that $E^\delta(\varepsilon) = \delta^2$, which completes the proof. $\square$

**Remark 2.1.** *In the particular case when the operator A has a dense range, that is* $(\text{Range } A)^\perp = \{0\}$*, we have of course* $g_\perp^\delta = 0$*, so that the assumption* $\|g_\perp^\delta\|_H < \delta$ *is automatically satisfied. This is why Theorem 2.1 is a generalization of the result given in [20] (see theorem 2.17).*

The context of Theorem 2.1 could be the following: we consider an idealized exact problem, which consists in, for the exact data $g \in H$, finding the exact solution $u \in V$ such that

$$Au = g.$$

However, in practice such exact problem is impossible to solve, because the exact data $g \in H$ is unknown. Instead, we measure some noisy data $g^\delta$ such that

$$\|g^\delta - g\|_H = \delta, \tag{7}$$

where $\delta$ is the amplitude of noise. The objective is to find a solution in $V$ from data $g^\delta$ which is close to $u$, in particular when such noisy data $g^\delta$ does not belong to the range of $A$. The Tikhonov regularization, which consists in computing $u_\varepsilon^\delta$ with the help of (2) for some $\varepsilon > 0$, is a classical way to approximate $u$. It is well-known that for $\delta = 0$, that is for exact data $g$, the corresponding Tikhonov solution $u_\varepsilon$ converges to $u$ in $V$ when $\varepsilon$ tends to 0. When $\delta > 0$ however, choosing $\varepsilon$ is not easy. In particular, each time $g^\delta \notin \text{Range } A$, the norm $\|u_\varepsilon^\delta\|_V$ tends to $+\infty$ when $\varepsilon$ tends to 0. The Morozov's principle is a classical way of choosing $\varepsilon$ such that (3) is satisfied. The general idea of the Morozov's principle is the following: since from (7) the data $g^\delta$ is corrupted by some noise of amplitude $\delta$, hence a solution to the problem $Au^\delta = g^\delta$ might fail to exist, it is not worth computing a Tikhonov solution such that $Au_\varepsilon^\delta \simeq g^\delta$ be satisfied with a better accuracy than $\delta$. More precisely, the Morozov's value $\varepsilon$ is chosen so that the error $\|Au_\varepsilon^\delta - g^\delta\|_H$ made in the resolution of the problem exactly coincides with the amplitude of noise $\delta$. A justification of the Morozov's rule is given by the following result. The proof is omitted since it is exactly the one given in theorem 2.17 of [20] in the restricted case when $A$ is a compact operator with dense range (in fact, the proof does not use these two assumptions).

**Proposition 2.1.** *With the same assumptions as in Theorem 2.1 and assumption (7), let us denote* $\varepsilon(\delta)$ *the value of* $\varepsilon$ *given by (3) and* $u_{\varepsilon(\delta)}^\delta$ *the corresponding solution to (2). Then*

$$\lim_{\delta \to 0} u_{\varepsilon(\delta)}^\delta = u \quad \text{in} \quad V.$$

**Remark 2.2.** *It should be noted that in the context which is described above, that is* $\|g^\delta - g\|_H = \delta$*, where* $g$ *is the exact data, the assumption* $\|g_\perp^\delta\|_H < \delta$ *is very unlikely to be violated. Indeed, let us denote* $r^\delta = g^\delta - g$ *the perturbation between the noisy and the exact data, and let us use the decomposition* $r^\delta = r_\parallel^\delta + r_\perp^\delta$*, with* $r_\parallel^\delta \in \overline{\text{Range } A}$ *and* $r_\perp^\delta \in (\text{Range } A)^\perp$*. Since* $g \in \text{Range } A$*, we have*

$$\|g_\perp^\delta\|_H^2 = \|r_\perp^\delta\|_H^2 = \delta^2 - \|r_\parallel^\delta\|_H^2,$$

*so that* $\|g_\perp^\delta\|_H < \delta$ *unless* $r_\parallel^\delta = 0$*, that is if and only if the perturbation* $r^\delta$ *only has a contribution in* $(\text{Range } A)^\perp$*.*

## 2.2 Interpretation of the Morozov's principle with duality in optimization

Now let us introduce and study a minimization problem which will be later on related to the Morozov's principle. We consider an orthogonal projector $P : H \to H$ on a closed subspace of $H$ such that $P$ is compact and $\operatorname{Range} P \subset \overline{\operatorname{Range} A}$. The operator $P$ will enable us to enforce the approximate solution to satisfy some *a priori* constraints, as will be clarified later.

We consider the problem

$$(\mathscr{P}_P^*) \quad \inf_{q \in H} G_P^\delta(q), \tag{8}$$

the functional $G_P^\delta$ being defined, for all $q \in H$, by

$$G_P^\delta(q) = \frac{1}{2}\|A^*q\|_V^2 + \delta\,\|(I-P)q\|_H - (g^\delta, q)_H, \tag{9}$$

where $I : H \to H$ is the identity operator and $A^* : H \to V$ is the adjoint operator of $A$.

**Lemma 2.1.** *The function $G_P^\delta$ is coercive if and only if the data $g^\delta$ satisfies $\|g_\perp^\delta\|_H < \delta$.*

*Proof.* Let us assume that $\|g_\perp^\delta\|_H < \delta$ and let us prove that $G_P^\delta(q) \to +\infty$ when $\|q\|_H \to +\infty$. Assume on the contrary that there exists some constant $C > 0$ and a sequence $(q_n)_{n \in \mathbb{N}}$ of elements in $H$ such that $\alpha_n = \|q_n\|_H \to +\infty$ while $G(q_n) \leq C$. Let us define $z_n = q_n/\|q_n\|_H$. Since $\|z_n\|_H = 1$, there exists a subsequence of $(z_n)_n$, still denoted $(z_n)_n$, such that $z_n \rightharpoonup z$ in $H$. We have

$$\frac{\alpha_n^2}{2}\|A^*z_n\|_V^2 + \alpha_n \delta \|(I-P)z_n\|_H - \alpha_n (g^\delta, z_n)_H \leq C, \tag{10}$$

hence

$$\frac{1}{2}\|A^*z_n\|_V^2 \leq \frac{1}{\alpha_n}(g^\delta, z_n)_H + \frac{C}{\alpha_n^2},$$

which implies that $A^*z_n \to 0$ in $V$ and since $A^*z_n \rightharpoonup A^*z$ in $V$, we get that $A^*z = 0$, that is $z \in \operatorname{Ker} A^* = (\operatorname{Range} A)^\perp$.

Another consequence of (10) is

$$\begin{aligned}
\delta\|(I-P)z_n\|_H &\leq (g^\delta, z_n)_H + \frac{C}{\alpha_n} \\
&= (g^\delta, z)_H - (g^\delta, z - z_n)_H + \frac{C}{\alpha_n} \\
&= (g_\perp^\delta, z)_H - (g^\delta, z - z_n)_H + \frac{C}{\alpha_n} \\
&= (g_\perp^\delta, (I-P)z)_H - (g^\delta, z - z_n)_H + \frac{C}{\alpha_n} \\
&\leq \|g_\perp^\delta\|_H \|(I-P)z\|_H + |(g^\delta, z - z_n)_H| + \frac{C}{\alpha_n},
\end{aligned}$$

where the second equality comes from the fact that $z \in (\operatorname{Range} A)^\perp$ and the third one from the fact that $\operatorname{Range} P \subset \overline{\operatorname{Range} A}$. Now, since $(I-P)z_n \rightharpoonup (I-P)z$ in $H$, we have

$$\begin{aligned}
\delta\|(I-P)z\|_H &\leq \liminf_{n \to +\infty} \delta\|(I-P)z_n\|_H \\
&\leq \liminf_{n \to +\infty} \left( \|g_\perp^\delta\|_H \|(I-P)z\|_H + |(g^\delta, z - z_n)_H| + \frac{C}{\alpha_n} \right) \\
&= \|g_\perp^\delta\|_H \|(I-P)z\|_H.
\end{aligned}$$

Since $\|g_\perp^\delta\|_H < \delta$, we obtain that $(I-P)z = 0$, in particular $z \in \operatorname{Range} P \subset \overline{\operatorname{Range} A}$, which together with $z \in (\operatorname{Range} A)^\perp$ implies $z = 0$.

Suppose finally that $(I-P)z_n \to 0$ in $H$. That $P$ is a compact operator and $z_n \rightharpoonup 0$ in $H$ yields $Pz_n \to 0$ in $H$, hence $z_n \to 0$ in $H$. But this contradicts the fact that $\|z_n\|_H = 1$ for all $n$. Then we can find a real $\varepsilon > 0$ and a subsequence of $(z_n)_n$, still denoted $(z_n)_n$, such that $\|(I-P)z_n\|_H \geq \varepsilon$ for all $n$. As a result, in view of (10),

$$C \geq \alpha_n(\delta\varepsilon - (g^\delta, z_n)_H) \to +\infty.$$

Such contradiction proves that $G_P^\delta$ is coercive if $\|g_\perp^\delta\|_H < \delta$.

Let us prove the converse statement, that is $\|g_\perp^\delta\|_H \geq \delta$, implies that $G_P^\delta$ is not coercive. For $\alpha \in \mathbb{R}$, by setting $q = \alpha g_\perp^\delta \in H$ in $G_P^\delta(q)$ we get

$$G_P^\delta(\alpha g_\perp) = \frac{\alpha^2}{2}\|A^* g_\perp^\delta\|_V + \alpha\delta\|(I-P)g_\perp^\delta\|_H - \alpha\|g_\perp^\delta\|_H^2.$$

We observe that $g_\perp^\delta \in (\operatorname{Range} A)^\perp = \operatorname{Ker} A^*$. The operator $P$ is self-adjoint as an orthogonal projector. Since in addition $\operatorname{Range} P \subset \overline{\operatorname{Range} A}$, we also observe that $g_\perp^\delta \in (\operatorname{Range} P)^\perp = \operatorname{Ker} P$. We conclude that

$$G_P^\delta(\alpha g_\perp) = \alpha(\delta\|g_\perp^\delta\|_H - \|g_\perp^\delta\|_H^2).$$

If $\|g_\perp^\delta\|_H > \delta$, then $\|g_\perp^\delta\|_H^2 > \delta\|g_\perp^\delta\|_H$, so that $G_P^\delta(\alpha g_\perp) \to -\infty$ when $\alpha \to +\infty$. If $\|g_\perp^\delta\|_H = \delta$, then $G_P^\delta(\alpha g_\perp) = 0$. In both cases, the functional $G_P^\delta$ is not coercive. □

From Lemma 2.1, by introducing the set

$$K_P^\delta = \{v \in V, \|Av - g^\delta\|_H \leq \delta,\ PAv = Pg^\delta\},$$

we obtain the following theorem.

**Theorem 2.2.** *If $\|g_\perp^\delta\|_H < \delta$, the optimization problem (8) has at least one solution $p^\delta \in H$. Then $u^\delta = A^* p^\delta$ belongs to the set $K_p^\delta$ and we have the identity*

$$\|u^\delta\|_V^2 + 2\, G_P^\delta(p^\delta) = 0. \tag{11}$$

*Proof.* The functional $G_P^\delta$ is continuous and convex on $H$. Existence of a minimizer $p^\delta$ of $G_P^\delta$ is then a consequence of the coercivity of $G_P^\delta$ by Lemma 2.1. Let us define $u^\delta = A^* p^\delta$. We consider two cases.

Let us firstly assume that $(I - P)p^\delta \neq 0$. The optimality of $G_P^\delta$ at $p^\delta$ writes

$$(A^*p^\delta, A^*q)_V + \frac{\delta}{\|(I-P)p^\delta\|_H}((I-P)p^\delta, (I-P)q)_H - (g^\delta, q)_H = 0, \quad \forall q \in H,$$

which implies, since $(I-P)^2 = I - P$,

$$Au^\delta - g^\delta = -\frac{\delta}{\|(I-P)p^\delta\|_H}(I-P)p^\delta. \tag{12}$$

We conclude that $\|Au^\delta - g^\delta\|_H = \delta$ and $PAu^\delta = Pg^\delta$, that is $u^\delta \in K_P^\delta$. In addition, we have

$$\begin{aligned}
\|u^\delta\|_V^2 + 2\, G_P^\delta(p^\delta) &= \|u^\delta\|_V^2 + \|A^*p^\delta\|_V^2 + 2\delta\|(I-P)p^\delta\|_H - 2(g^\delta, p^\delta)_H \\
&= 2(Au^\delta - g^\delta, p^\delta)_H + 2\delta\|(I-P)p^\delta\|_H \\
&= -\frac{2\delta}{\|(I-P)p^\delta\|_H}((I-P)p^\delta, p^\delta)_H + 2\delta\|(I-P)p^\delta\|_H \\
&= 0,
\end{aligned}$$

where the third equality is a consequence of (12).

Let us on the contrary assume that $(I - P)p^\delta = 0$. Then the functional $G_P^\delta$ is not differentiable at point $p^\delta$. The optimality however writes $0 \in \partial G_P^\delta(p^\delta)$, where $\partial G_P^\delta(p^\delta)$ denotes the subdifferential of $G_P^\delta$ at point $p^\delta$. In view of (9), and by using the classical rules for subdifferential computations, we have

$$\partial G_P^\delta(p^\delta) = AA^*p^\delta - g^\delta + \delta(I-P)\partial(\|\cdot\|_H)(0).$$

Since $\partial(\|\cdot\|_H)(0)$ is the unit ball of $H$ centered at 0, that $0 \in \partial G_P^\delta(p^\delta)$ implies that $\|Au^\delta - g^\delta\|_H \leq \delta$. In addition, we have

$$G_P^\delta(p^\delta) = \inf_{q \in H} G_P^\delta(q) \leq \inf_{q \in \operatorname{Range} P} G_P^\delta(q) = \inf_{q \in \operatorname{Range} P} \tilde{G}^\delta(q),$$

where

$$\tilde{G}^\delta(q) = \frac{1}{2}\|A^*q\|_V^2 - (g^\delta, q)_H.$$

Since $p^\delta \in \operatorname{Range} P$, in fact we have
$$G_P^\delta(p^\delta) = \tilde{G}^\delta(p^\delta) = \inf_{q \in \operatorname{Range} P} \tilde{G}^\delta(q).$$

The optimality of $\tilde{G}^\delta$ in $\operatorname{Range} P$ at point $p^\delta$ writes
$$(A^*p^\delta, A^*q)_V - (g^\delta, q)_H = 0, \quad \forall q \in \operatorname{Range} P,$$

which amounts to
$$Au^\delta - g^\delta \in (\operatorname{Range} P)^\perp = \operatorname{Ker} P. \tag{13}$$

We conclude that $u^\delta \in K_P^\delta$. In addition,
$$\begin{aligned} \|u^\delta\|_V^2 + 2\,G_P^\delta(p^\delta) &= \|u^\delta\|_V^2 + \|A^*p^\delta\|_V^2 - 2(g^\delta, p^\delta)_H \\ &= 2(Au^\delta - g^\delta, p^\delta)_H \\ &= 0, \end{aligned}$$

the last equality being a consequence of (13). The proof is complete. □

In order to give a precise meaning to $u^\delta$, we need the following lemma.

**Lemma 2.2.** *For all $v \in K_P^\delta \setminus \{u^\delta\}$, it holds that $\|v\|_V > \|u^\delta\|_V$.*

*Proof.* Assume that $v \in K_P^\delta$ and let us introduce $q = g^\delta - Av$, which satisfies $\|q\|_H \leq \delta$ and $Pq = 0$. We have
$$\begin{aligned} \frac{1}{2}(\|v\|_V^2 - \|u^\delta\|_V^2) &= \frac{1}{2}\|v\|_V^2 + G_P^\delta(p^\delta) \\ &= \frac{1}{2}\|v\|_V^2 + \frac{1}{2}\|A^*p^\delta\|_V^2 + \delta\|(I-P)p^\delta\|_H - (g^\delta, p^\delta)_H \\ &= \frac{1}{2}\|v\|_V^2 + \frac{1}{2}\|u^\delta\|_V^2 + \delta\|(I-P)p^\delta\|_H - (Av, p^\delta)_H - (q, p^\delta)_H \\ &= \frac{1}{2}\|v\|_V^2 + \frac{1}{2}\|u^\delta\|_V^2 - (v, u^\delta)_V + \delta\|(I-P)p^\delta\|_H - ((I-P)q, p^\delta)_H \\ &= \frac{1}{2}\|v - u^\delta\|_V^2 + \delta\|(I-P)p^\delta\|_H - (q, (I-P)p^\delta)_H \\ &\geq \frac{1}{2}\|v - u^\delta\|_V^2, \end{aligned}$$

where the first equality is a consequence of (11), while the fourth one is a consequence of $Pq = 0$ and the last inequality uses that $\|q\|_H \leq \delta$. The proof is complete. □

**Theorem 2.3.** *If $\|g_\perp^\delta\|_H < \delta$, for any solution $p^\delta \in H$ to the minimization problem (8), $u^\delta = A^*p^\delta \in V$ is the unique solution to the minimization problem*
$$(\mathscr{P}_P) \quad \inf_{v \in K_P^\delta} \|v\|_V^2. \tag{14}$$

*Proof.* We first remark that the optimization problem (14) has a unique solution. Indeed, the set $K_P^\delta$ is convex, closed and non-empty (it contains $u^\delta$ from Theorem 2.2). In addition, the cost function to minimize is coercive and strictly convex. Lemma 2.2 shows that $u^\delta$ coincides with the unique solution to such problem (14). □

**Remark 2.3.** *An important consequence of Theorem 2.3 is that $u^\delta = A^*p^\delta$ does not depend on the solution $p^\delta$ to the optimization problem (8).*

The following result specifies in which sense $u^\delta$ is an approximation of the exact solution $u$.

**Theorem 2.4.** *Let us assume that $g^\delta$ satisfies (1) and (7). If in addition $Pg^\delta = Pg$, then*
$$\lim_{\delta \to 0} u^\delta = u \quad \text{in} \quad V.$$

*Proof.* We have
$$\|Au - g^\delta\|_H = \|g - g^\delta\|_H = \delta$$
and $PAu = Pg^\delta$, so that $u \in K_P^\delta$. By Lemma 2.2 we obtain that $\|u^\delta\|_V \leq \|u\|_V$. From the sequence $(u^\delta)_\delta$ in $V$, we can then extract a subsequence still denoted $(u^\delta)_\delta$, such that $u^\delta \rightharpoonup w$ in $V$ when $\delta \to 0$. Besides, we have both $Au^\delta \to g$ in $H$ and $Au^\delta \rightharpoonup Aw$ in $H$, hence $Aw = g$, that is $w = u$ from the injectivity of $A$. As a result, $u^\delta \rightharpoonup u$ in $V$. Then
$$\|u^\delta - u\|_V^2 = \|u^\delta\|_V^2 + \|u\|_V^2 - 2(u^\delta, u)_V \leq 2(u - u^\delta, u)_V,$$
and we get that $u^\delta \to u$ in $V$. We easily conclude that all the sequence $(u^\delta)_\delta$, and not only a subsequence, converges to $u$ in $V$. □

We wish now, in the particular case when $P = 0$, to relate the solutions $p^\delta$ which minimize the functional $G_0^\delta$ (note that $G_0^\delta$ coincides with the functional $G_P^\delta$ given by (9) when $P = 0$) to the classical solution of the Tikhonov regularized problem associated with operator $A$ when the regularization parameter is chosen according to the Morozov's discrepancy principle. We hence consider the minimization problem

$$(\mathscr{P}^*) \quad \inf_{q \in H} G_0^\delta(q) = \inf_{q \in H} \left( \frac{1}{2} \|A^*q\|_V^2 + \delta \|q\|_H - (g^\delta, q)_H \right). \tag{15}$$

The solutions to problem $(\mathscr{P}^*)$ will enable us to obtain a practical method to compute the Morozov's value $\varepsilon(\delta) > 0$ given by Theorem 2.1 and the corresponding Tikhonov solution $u_\varepsilon^\delta$. In the appendix, we show how the problem (15) can be derived in a constructive way by using the theory exposed in [17].

**Theorem 2.5.** *If the noisy data $g^\delta$ satisfies the assumption (1), the problem $(\mathscr{P}^*)$ given by (15) has at least one solution, and $u^\delta = A^*p^\delta$ coincides with the Tikhonov solution $u_\varepsilon^\delta$ of problem (2) where $\varepsilon(\delta)$ is the unique value of $\varepsilon > 0$ such that $\|Au_\varepsilon^\delta - g^\delta\|_H = \delta$ according to Theorem 2.1. Lastly, $p^\delta \neq 0$ and*
$$\varepsilon(\delta) = \frac{\delta}{\|p^\delta\|_H}.$$

*Proof.* From Theorem 2.2, we already know that the problem $(\mathscr{P}^*)$ has solutions which are denoted $p^\delta$. Let us verify that $p^\delta \neq 0$. Actually, let us take $g_1^\delta = g^\delta / \|g^\delta\|_H$. For $\varepsilon > 0$, we have
$$G_0^\delta(\varepsilon g_1^\delta) = \frac{\varepsilon^2}{2} \|A^* g_1^\delta\|_V^2 + \varepsilon(\delta - \|g^\delta\|_H).$$

If $\varepsilon$ is sufficiently small, $G_0^\delta(\varepsilon g_1^\delta)$ has the sign of $\delta - \|g^\delta\|_H < 0$, hence there exists $q \in H$ such that $G_0^\delta(q) < 0 = G_0^\delta(0)$, and the solutions $p^\delta$ do not vanish. From the proof of Theorem 2.2 in the case $P = 0$, in particular in view of (12), if we denote $u^\delta = A^* p^\delta$, we directly obtain that $u^\delta$ satisfies the Morozov's principle
$$\|Au^\delta - g^\delta\|_H = \delta$$
and the equation
$$A^*(Au^\delta) + \frac{\delta}{\|p^\delta\|_H} u^\delta = A^* g^\delta.$$

We hence conclude that if we take $\varepsilon(\delta) = \delta / \|p^\delta\|_H$, the function $u^\delta$ satisfies $A^*(Au^\delta) + \varepsilon(\delta) u^\delta = A^* g^\delta$, which means that it is the unique solution $u_\varepsilon^\delta$ to the problem (2) associated with that $\varepsilon(\delta)$ and furthermore satisfies the Morozov's principle. □

**Remark 2.4.** *In the case when $P = 0$, the above theorem provides a strategy to find the Tikhonov/Morozov solution associated with noisy data $g^\delta$. It consists in finding first a solution to problem $(\mathscr{P}^*)$. The Morozov solution is then obtained by applying $A^*$ to any solution of $(\mathscr{P}^*)$.*

*In the case when the operator $P$ is not 0, for $p^\delta$ the solutions to problem $(\mathscr{P}_P^*)$, the corresponding solutions $u^\delta = A^* p^\delta$ can not be related to the Tikhonov problem (2). However, since $\|Au^\delta - g^\delta\|_H \leq \delta$, they satisfy the Morozov's principle in the sense of an inequality instead of an equality when $\delta = \|g^\delta - g\|_H$. The role of the operator $P$ is to impose that $PAu^\delta = Pg^\delta$, which ensures that some particular reliable features of the data $g^\delta$ are satisfied exactly by $Au^\delta$. For example, the noise often affects the high frequencies of the measurements. Hence we are tempted to be more confident in the low frequencies of the data than in their high frequencies. It is then natural to impose that a finite number of low frequency components of the data be exactly satisfied by the approximate solution, which can be achieved by using a specific operator $P$. We will present some examples of projector $P$ in the case of the data assimilation problem for the Laplace equation.*

# 3 An abstract framework for data assimilation problems

In this section we introduce a general framework for a class of data assimilation problems. The three applications that we will present in the next section are particular cases of such general framework. Let us consider $V$, $M$ and $O$ three Hilbert spaces, $b$ a bilinear continuous mapping on $V \times M$ and the corresponding operator $B : V \to M$ such that
$$(Bu, \lambda)_M = b(u, \lambda), \quad \forall (u, \lambda) \in V \times M,$$
as well as a continuous operator $C : V \to O$. We assume that the operator $A : V \to H = M \times O$ such that $Au = (Bu, Cu)$ is injective. We formulate our abstract data assimilation problem as follows: for data $f \in O$, find $u \in V$ such that $Bu = 0$ and $Cu = f$.

By the injectivity of $A$, such problem has at most one solution but in many situations, it is ill-posed because $A$ is not onto. This is why we propose, for $\varepsilon > 0$, the following regularized weak formulation: for $f \in O$, find $(u_\varepsilon, \lambda_\varepsilon) \in V \times M$ such that for all $(v, \mu) \in V \times M$,

$$\begin{cases} \varepsilon(u_\varepsilon, v)_V + (Cu_\varepsilon, Cv)_O + b(v, \lambda_\varepsilon) = (f, Cv)_O \\ b(u_\varepsilon, \mu) - (\lambda_\varepsilon, \mu)_M = 0. \end{cases} \quad (16)$$

As recalled in the introduction, the principal motivation for introducing such variational mixed formulation is to find an approximate solution to our ill-posed problem in the space $V$, which is the natural space of the true solution $u$. In addition, when it comes to the discretization with the Finite Element Method, it enables us to consider simple conforming finite elements. Alternatively, if we directly apply the ideas of [21], for instance, the approximate solution has to be searched in a space more regular than $V$ and the discretization requires some more cumbersome finite elements, as can be seen in [7].

**Remark 3.1.** *We have chosen here to restrict ourselves to the homogeneous equation $Bu = 0$ instead of a non-homogeneous equation $Bu = \ell$ for some $\ell \in M$, in order to handle one single data $f$ instead of a couple of data $(\ell, f)$, and hence simplify the analysis. The fully non-homogeneous case is for example addressed in [8], where a Cauchy problem for the Helmholtz equation in the presence of both noisy Dirichlet and Neumann data is considered. This increases the difficulty to apply the Mororov's principle since the two Cauchy data are independently perturbed by some noise while a single regularization parameter is at one's disposal.*

The weak formulation (16) is justified by the following theorem.

**Theorem 3.1.** *For all $\varepsilon > 0$, the weak formulation (16) has a unique solution. Furthermore, if there exists $u \in V$ such that $Bu = 0$ and $Cu = f$, then $(u_\varepsilon, \lambda_\varepsilon) \to (u, 0)$ in $V \times M$ when $\varepsilon \to 0$.*

*Proof.* The weak formulation (16) is equivalent to: find $(u_\varepsilon, \lambda_\varepsilon) \in V \times M$ such that for all $(v, \mu) \in V \times M$,
$$\mathcal{A}((u_\varepsilon, \lambda_\varepsilon); (v, \mu)) = \mathcal{L}((v, \mu)),$$
with
$$\mathcal{A}((u, \lambda); (v, \mu)) = \varepsilon(u, v)_V + (Cu, Cv)_O + b(v, \lambda) - b(u, \mu) + (\lambda, \mu)_M, \quad \mathcal{L}((v, \mu)) = (f, Cv)_O.$$
If suffices to apply the Lax-Milgram Lemma, the coercivity of $\mathcal{A}$ being ensured by
$$\mathcal{A}((u, \lambda); (u, \lambda)) = \varepsilon\|u\|_V^2 + \|Cu\|_O^2 + \|\lambda\|_M^2 \geq \min(\varepsilon, 1)(\|u\|_V^2 + \|\lambda\|_M^2).$$

Now assume that there exists $u \in V$ such that $Bu = 0$ and $Cu = f$. By the injectivity of $A$, such $u$ is uniquely defined. Since $Bu = 0$ and $Cu = f$, the system (16) implies that for all $(v, \mu) \in V \times M$,
$$\begin{cases} \varepsilon(u_\varepsilon, v)_V + (C(u_\varepsilon - u), Cv)_O + b(v, \lambda_\varepsilon) = 0 \\ b(u_\varepsilon - u, \mu) - (\lambda_\varepsilon, \mu)_M = 0. \end{cases}$$
Choosing $v = u_\varepsilon - u$ in the first equation and $\mu = \lambda_\varepsilon$ in the second equation, taking the difference of the two obtained equations implies that for all $\varepsilon > 0$,
$$\varepsilon(u_\varepsilon - u, u_\varepsilon)_V + \|C(u_\varepsilon - u)\|_O^2 + \|\lambda_\varepsilon\|_M^2 = 0. \quad (17)$$

Identity (17) implies that $(u_\varepsilon - u, u_\varepsilon)_V \leq 0$, hence $(u_\varepsilon)_\varepsilon$ is bounded in $V$. There exists a subsequence of $(u_\varepsilon)_\varepsilon$, still denoted $(u_\varepsilon)_\varepsilon$, which weakly converges to some $w \in V$. From (17), we also deduce that $(Cu_\varepsilon)_\varepsilon$ converges to $Cu$ in $O$ and that $(\lambda_\varepsilon)$ converges to 0 in $M$. From the second equation of (16), we have that $Bu_\varepsilon = \lambda_\varepsilon$, which implies that $(Bu_\varepsilon)_\varepsilon$ converges to 0. Since the sequences $(Bu_\varepsilon)_\varepsilon$ and $(Cu_\varepsilon)_\varepsilon$ weakly converges to $Bw$ in $M$ and $Cw$ in $O$, respectively, we obtain that $Bw = 0 = Bu$ and $Cw = Cu$. The injectivity of $A = (B,C)$ yields $w = u$. It remains to remark that

$$\|u_\varepsilon - u\|_V^2 = (u_\varepsilon - u, u_\varepsilon)_V - (u, u_\varepsilon - u)_V \leq -(u, u_\varepsilon - u)_V,$$

which implies that the weak convergence of $(u_\varepsilon)$ implies the strong convergence of $(u_\varepsilon)$ in $V$. A classical contradiction argument proves that all the sequence $(u_\varepsilon)_\varepsilon$, and not only the subsequence, tends to $u$ in $V$ when $\varepsilon$ tends to 0. $\square$

**Remark 3.2.** *Reading the above proof carefully shows that Theorem 3.1 also holds if the mapping $v \in V \mapsto \|v\|_V \in \mathbb{R}_+$ is no more a norm in $V$ but only a semi-norm, provided the mapping $v \in V \mapsto (\|v\|_V^2 + \|Cv\|_O^2)^{1/2} \in \mathbb{R}_+$ is a norm in $V$ such that $V$, equipped with such a norm, is complete.*

We now offer a link between the weak formulation (16) and the Tikhonov regularization (2). For $f \in O$, finding $u \in V$ such that $Bu = 0$ and $Cu = f$ is equivalent to finding $u \in V$ such that $Au = g$, where $g = (0, f) \in H$.

**Proposition 3.1.** *The pair $(u_\varepsilon, \lambda_\varepsilon) \in V \times M$ is the solution to problem (16) if and only if $u_\varepsilon$ is the solution to the problem: find $u_\varepsilon \in V$ such that for all $v \in V$,*

$$(Au_\varepsilon, Av)_H + \varepsilon(u_\varepsilon, v)_V = (g, Av)_H \tag{18}$$

*for $g = (0, f) \in H$, and $\lambda_\varepsilon = Bu_\varepsilon$.*

*Proof.* The function $u_\varepsilon \in V$ satisfies the problem (18) if and only if

$$(Bu_\varepsilon, Bv)_M + (Cu_\varepsilon, Cv)_O + \varepsilon(u_\varepsilon, v)_V = (0, Bv)_M + (f, Cv)_O, \quad \forall v \in V,$$

that is, setting $\lambda_\varepsilon = Bu_\varepsilon$,

$$\varepsilon(u_\varepsilon, v)_V + (Cu_\varepsilon, Cv)_O + (Bv, \lambda_\varepsilon)_M = (f, Cv)_O, \quad \forall v \in V,$$

and

$$(Bu_\varepsilon, \mu)_M - (\lambda_\varepsilon, \mu)_M = 0, \quad \forall \mu \in M,$$

which is equivalent to formulation (16). $\square$

That the problem $Au = g$ is ill-posed and the above equivalence between the Tikhonov regularization for the operator $A$ and the mixed formulation (16) prevents us from setting $\varepsilon = 0$ in (16). In particular, such ill-posedness implies that the bilinear form $b$ does not satisfy the inf-sup condition, as shown in [3] and [8]. If we now consider, for $\delta > 0$, some noisy data $f^\delta \in O$, let us denote $(u_\varepsilon^\delta, \lambda_\varepsilon^\delta) \in V \times M$ the solution to problem (16) associated with data $f^\delta$. From Proposition 3.1, we get that $u_\varepsilon^\delta$ is the solution to problem (18) for $g^\delta = (0, f^\delta)$ and that $\lambda_\varepsilon^\delta = Bu_\varepsilon^\delta$. From Theorem 2.1 applied to operator $A = (B, C)$ from $V$ to $H = M \times O$, we immediately obtain the following result.

**Corollary 3.1.** *Let us denote $(\lambda_\perp^\delta, f_\perp^\delta) \in M \times O$ the orthogonal projection of data $(0, f^\delta) \in M \times O$ on $(\text{Range } A)^\perp$ and assume that*

$$\|\lambda_\perp^\delta\|_M^2 + \|f_\perp^\delta\|_O^2 < \delta^2 < \|f^\delta\|_O^2.$$

*There exists a unique $\varepsilon > 0$ such that*

$$\|\lambda_\varepsilon^\delta\|_M^2 + \|Cu_\varepsilon^\delta - f^\delta\|_O^2 = \delta^2.$$

**Remark 3.3.** *Corollary 3.1 is close to Theorem 2.10 in [8]. However, here we point out that the statement of such theorem is not correct in the sense that the following assumption should have been added: using the notations of [8], the operator $\mathcal{A} = (A, B)$ shall have a dense range. In particular, that $\mathcal{A}$ has a dense range is not a consequence of the fact that the operators $A$ and $B$ both have a dense range, contrary to what is claimed in [8]. Fortunately, the property that the operator $\mathcal{A}$ has a dense range is actually true in the particular case considered in [8], which is very similar to the Cauchy problem for the Laplace equation addressed hereafter (see section 4.2, in particular Lemma 4.5).*

In order to compute the value of $\varepsilon$ and the Morozov's solution $u_\varepsilon^\delta$ given by Corollary 3.1 by using duality, we wish to give a more explicit form of the minimization problem (8) as well as the Fréchet derivative of $G_P^\delta$ in the data assimilation framework. This Fréchet derivative will be required to solve the minimization problem with the help of an iterative gradient method. More generally, we consider two orthogonal projectors $P_M : M \to M$ and $P_O : O \to O$ on a closed subspace of $M$ and $O$, respectively, such that $P_M$ and $P_O$ are compact. We hence obtain an orthogonal projector $P = (P_M, P_O) : M \times O \to M \times O$ which is compact and assume in addition that Range $P \subset \overline{\text{Range } A}$. We are hence in a position to state the following proposition.

**Proposition 3.2.** *Let us identify the spaces $V^*$, $M^*$, $O^*$ with $V$, $M$, $O$, respectively. The dual problem (8) reduces to*

$$(\mathscr{P}_P^*) \quad \inf_{(\lambda^*, f^*) \in M \times O} G_P^\delta(\lambda^*, f^*) \tag{19}$$

*with for $(\lambda^*, f^*) \in M \times O$,*

$$G_P^\delta(\lambda^*, f^*) = \frac{1}{2}\|B^*\lambda^* + C^*f^*\|_V^2 + \delta\,\|(I-P)(\lambda^*, f^*)\|_{M \times O} - (f^\delta, f^*)_O. \tag{20}$$

*Let us define $u^* \in V$ by the weak formulation*

$$(u^*, v)_V = b(v, \lambda^*) + (f^*, Cv)_O, \quad \forall v \in V. \tag{21}$$

*For $(I - P)(\lambda^*, f^*) \neq 0$, the two partial Fréchet derivatives $\partial_\lambda G_P^\delta \in M$ and $\partial_f G_P^\delta \in O$ of $G_P^\delta$ at point $(\lambda^*, f^*)$ are given by the weak formulation*

$$(\partial_\lambda G_P^\delta, \mu)_M = b(u^*, \mu) + \frac{\delta}{\|(I-P)(\lambda^*, f^*)\|_{M \times O}}((I_M - P_M)\lambda^*, \mu)_M, \quad \forall \mu \in M \tag{22}$$

*and the identity*

$$\partial_f G_P^\delta = Cu^* + \frac{\delta}{\|(I-P)(\lambda^*, f^*)\|_{M \times O}}(I_O - P_O)f^* - f^\delta, \tag{23}$$

*respectively.*

*Proof.* Formula (20) is obtained having in mind that $g^\delta = (0, f^\delta)$, for $f^\delta \in O$, and by observing that since $Av = (Bv, Cv)_{M \times O}$ for all $v \in V$, we have for all $(\lambda^*, f^*) \in M^* \times O^*$, $A^*(\lambda^*, f^*) = B^*\lambda^* + C^*f^*$. In addition, $(I - P)(\lambda^*, f^*)) = ((I_M - P_M)\lambda^*, (I_O - P_O)f^*)$. Differentiating (20), the partial derivative $\partial_\lambda G_P^\delta$ is given, introducing $u^* = B^*\lambda^* + C^*f^* \in V$, by

$$(\partial_\lambda G_P^\delta, \mu)_M = (u^*, B^*\mu)_V + \frac{\delta}{\|(I-P)(\lambda^*, f^*)\|_{M \times O}}((I_M - P_M)\lambda^*, \mu)_M,$$

which implies (21) and (22), while the partial derivative $\partial_f G_P^\delta$ is given, for $h \in O$, by

$$(\partial_f G_P^\delta, h)_O = (u^*, C^*h)_V + \frac{\delta}{\|(I-P)(\lambda^*, f^*)\|_{M \times O}}((I_O - P_O)f^*, h)_O - (f^\delta, h)_O,$$

which implies (23). □

## 4 Some applications

### 4.1 The data assimilation problem for the Laplace equation

Let us consider a bounded Lipschitz domain $\Omega \subset \mathbb{R}^d$, $d > 1$, and a subdomain $\omega \Subset \Omega$ of class $C^2$. The data assimilation problem for the Laplace equation consists here, for some data $f \in L^2(\omega)$, in finding $u \in H^1(\Omega)$ such that $\Delta u = 0$ in $\Omega$ (in the sense of distributions) and $u|_\omega = f$. Clearly, such problem has at most one solution in virtue of Holmgren's theorem, but is ill-posed, since any harmonic solution in $\Omega$ is infinitely smooth in $\Omega$.

The data assimilation problem for the Laplace equation is a particular case of the general framework exposed in section 3. It corresponds to spaces $V = H^1(\Omega)$, $M = H_0^1(\Omega)$, $O = L^2(\omega)$, the bilinear form $b$ on $H^1(\Omega) \times H_0^1(\Omega)$ such that $b(u, \lambda) = \int_\Omega \nabla u \cdot \nabla \lambda\, dx$, and the associated operator $B : H^1(\Omega) \to H_0^1(\Omega)$ such that $(Bu, \lambda)_M = b(u, \lambda)$ for all $(u, \lambda) \in H^1(\Omega) \times H_0^1(\Omega)$, while the operator $C : H^1(\Omega) \to L^2(\omega)$ is the restriction operator. The space $H_0^1(\Omega)$ is equipped with the $H^1$ semi-norm.

**Lemma 4.1.** *For $f \in L^2(\omega)$, $u \in H^1(\Omega)$ satisfies the data assimilation problem for the Laplace equation if and only if $Bu = 0$ and $Cu = f$.*

*Proof.* The solution $u \in H^1(\Omega)$ satisfies the data assimilation problem iff $u|_\omega = f$ and
$$b(u, \mu) = \int_\Omega \nabla u \cdot \nabla \mu \, dx = 0, \quad \forall \mu \in H_0^1(\Omega),$$
in other words iff $Bu = 0$ and $Cu = f$. □

The injectivity of the operator $A = (B, C)$ is a direct consequence of Holmgren's Theorem. The corresponding formulation (16) is: for $f \in L^2(\omega)$, find $(u_\varepsilon, \lambda_\varepsilon) \in H^1(\Omega) \times H_0^1(\Omega)$ such that for all $(v, \mu) \in H^1(\Omega) \times H_0^1(\Omega)$,
$$\begin{cases} \varepsilon \int_\Omega (u_\varepsilon v + \nabla u_\varepsilon \cdot \nabla v) \, dx + \int_\omega u_\varepsilon v \, dx + \int_\Omega \nabla v \cdot \nabla \lambda_\varepsilon \, dx = \int_\omega f v \, dx, \\ \int_\Omega \nabla u_\varepsilon \cdot \nabla \mu \, dx - \int_\Omega \nabla \lambda_\varepsilon \cdot \nabla \mu \, dx = 0. \end{cases} \quad (24)$$

A Poincaré-type inequality implies that the mapping $v \mapsto (\int_\Omega |\nabla v|^2 \, dx + \int_\omega v^2 \, dx)^{1/2}$ is an equivalent norm to the standard norm on $H^1(\Omega)$. In view of Remark 3.2, the term $\int_\Omega u_\varepsilon v \, dx$ in (24) could therefore be dropped. With a view to applying Corollary 3.1, let us identify $(\text{Range } A)^\perp$, which is the aim of the next Lemma. In what follows, we identify the space $H_0^2(\omega)$ with the space of functions in $H_0^2(\Omega)$ which vanish in $\Omega \setminus \overline{\omega}$.

**Lemma 4.2.** *We have that*
$$(\text{Range } A)^\perp = \{(\mu, h) \in H_0^2(\omega) \times L^2(\omega), \quad \Delta \mu = h\}.$$

*Proof.* Let us consider $(\mu, h) \in H_0^1(\Omega) \times L^2(\omega)$. We have $(\mu, h) \in (\text{Range } A)^\perp$ if and only if for all $v \in H^1(\Omega)$
$$(Bv, \mu)_M + (Cv, h)_O = 0,$$
that is
$$\int_\Omega \nabla v \cdot \nabla \mu \, dx + \int_\omega v h \, dx = 0,$$
which is equivalent to
$$\begin{cases} \Delta \mu = 1_\omega h & \text{in } \Omega \\ \partial_\nu \mu = 0 & \text{on } \partial \Omega. \end{cases} \quad (25)$$
Here, $1_\omega$ is the indicator function of the set $\omega$. Hence, $(\mu, h) \in (\text{Range } A)^\perp$ implies that $\mu$ satisfies $\Delta \mu = 0$ in $\Omega \setminus \overline{\omega}$, $\mu = 0$ and $\partial_\nu \mu = 0$ on $\partial \Omega$, hence $\mu = 0$ in $\Omega \setminus \overline{\omega}$ from uniqueness of the Cauchy problem. We obtain that $\mu = 0$ and $\partial_\nu \mu = 0$ on $\partial \omega$, and from a standard regularity result and the fact that $\omega$ is of class $C^2$, we conclude that $\mu \in H_0^2(\omega)$. Conversely, it is straightforward to check that if $(\mu, h) \in \{(\mu, h) \in H_0^2(\omega) \times L^2(\omega), \Delta \mu = h\}$, then in particular $(\mu, h) \in H_0^1(\Omega) \times L^2(\omega)$ and satisfies (25), which finally yields $(\mu, h) \in (\text{Range } A)^\perp$. □

In the following lemma, we specify the orthogonal projection of $(0, f) \in H_0^1(\Omega) \times L^2(\omega)$ on $(\text{Range } A)^\perp$.

**Lemma 4.3.** *The projection of $(0, f) \in H_0^1(\Omega) \times L^2(\omega)$ on $(\text{Range } A)^\perp$ is the pair $(\lambda_\perp, f_\perp) \in H_0^1(\Omega) \times L^2(\omega)$, where $\lambda_\perp$ is the unique solution in $H_0^2(\omega)$ of the weak formulation set in the subdomain $\omega$:*
$$\int_\omega \Delta \lambda_\perp \, \Delta \mu \, dx + \int_\omega \nabla \lambda_\perp \cdot \nabla \mu \, dx = \int_\omega f \, \Delta \mu \, dx, \quad \forall \mu \in H_0^2(\omega), \quad (26)$$
*and $f_\perp = \Delta \lambda_\perp$.*

*Proof.* Let us find the orthogonal projection $(\lambda_\perp, f_\perp) \in H_0^1(\Omega) \times L^2(\omega)$ of any pair $(\lambda, f) \in H_0^1(\Omega) \times L^2(\omega)$ on $(\text{Range } A)^\perp$. Such orthogonal projection is characterized by
$$\begin{cases} (\lambda_\perp, f_\perp) \in (\text{Range } A)^\perp, \\ (\lambda - \lambda_\perp, f - f_\perp) \perp (\mu, h), \quad \forall (\mu, h) \in (\text{Range } A)^\perp, \end{cases}$$

which from Lemma 4.2 is equivalent to

$$\begin{cases} \lambda_\perp \in H_0^2(\omega) \quad \text{and} \quad f_\perp = \Delta \lambda_\perp, \\ \int_\Omega \nabla(\lambda - \lambda_\perp) \cdot \nabla\mu\, dx + \int_\omega (f - f_\perp) h\, dx = 0, \quad \forall \mu \in H_0^2(\omega) \quad \text{and} \quad h = \Delta\mu, \end{cases}$$

and then to

$$\begin{cases} \lambda_\perp \in H_0^2(\omega) \quad \text{and} \quad f_\perp = \Delta \lambda_\perp, \\ \int_\omega \Delta\lambda_\perp \Delta\mu\, dx + \int_\omega \nabla\lambda_\perp \cdot \nabla\mu\, dx = \int_\omega f \Delta\mu\, dx + \int_\omega \nabla\lambda \cdot \nabla\mu\, dx, \quad \forall \mu \in H_0^2(\omega). \end{cases} \quad (27)$$

The result follows by setting $\lambda = 0$. $\square$

**Remark 4.1.** *It should be noted that (26) is a fourth-order problem, which from the numerical point of view requires some cumbersome finite elements. In the case of the heat equation, computing the corresponding projection on* $(\mathrm{Range}\, A)^\perp$ *is even more complicated, as can be seen on (36).*

We can now apply Corollary 3.1 to the data assimilation problem. Let $f^\delta \in L^2(\omega)$ and $(\lambda_\perp^\delta, f_\perp^\delta)$ the orthogonal projection of $(0, f^\delta)$ on $(\mathrm{Range}\, A)^\perp$, which can be computed with the help of Lemma 4.3. If we assume that

$$\|\lambda_\perp^\delta\|_{H_0^1(\Omega)}^2 + \|f_\perp^\delta\|_{L^2(\omega)}^2 < \delta^2 < \|f^\delta\|_{L^2(\omega)}^2,$$

then there exists a unique $\varepsilon > 0$ such that

$$\int_\Omega |\nabla \lambda_\varepsilon^\delta|^2\, dx + \int_\omega |u_\varepsilon^\delta - f^\delta|^2\, dx = \delta^2,$$

where $(u_\varepsilon^\delta, \lambda_\varepsilon^\delta)$ is the unique solution to the problem (24) for data $f^\delta$.

**Remark 4.2.** *By choosing $\mu = \lambda_\perp$ in the weak formulation (26), we notice that*

$$\|\lambda_\perp^\delta\|_{H_0^1(\Omega)}^2 + \|f_\perp^\delta\|_{L^2(\omega)}^2 = \int_\omega f^\delta f_\perp^\delta\, dx.$$

In order to compute the Morozov parameter $\varepsilon$ and the corresponding solution $u_\varepsilon^\delta$ which are given above, one may be tempted to solve the minimization problem given by (19) (20) following the strategy exposed in Remark 2.4. Let us see how Proposition 3.2 is specified in the case of the data assimilation problem for the Laplace equation, in particular as concerns the computation of the Fréchet derivative of $G_P^\delta$. In view of Proposition 3.2, the solution $u^* \in H^1(\Omega)$ is given, for $(\lambda^*, f^*) \in H_0^1(\Omega) \times L^2(\omega)$, by the weak formulation

$$\int_\Omega (u^* v + \nabla u^* \cdot \nabla v)\, dx = \int_\Omega \nabla v \cdot \nabla \lambda^*\, dx + \int_\omega f^* v\, dx, \quad \forall v \in H^1(\Omega). \quad (28)$$

The dual problem (8) then amounts to

$$(\mathscr{P}_P^*) \quad \inf_{(\lambda^*, f^*) \in H_0^1(\Omega) \times L^2(\omega)} G_P^\delta(\lambda^*, f^*) \quad (29)$$

with for $(\lambda^*, f^*) \in H_0^1(\Omega) \times L^2(\omega)$,

$$G_P^\delta(\lambda^*, f^*) = \frac{1}{2}\|u^*(\lambda^*, f^*)\|_{H^1(\Omega)}^2 + \delta\,\|(I - P)(\lambda^*, f^*)\|_{H_0^1(\Omega) \times L^2(\omega)} - \int_\omega f^\delta f^*\, dx, \quad (30)$$

with

$$\|(I - P)(\lambda^*, f^*)\|_{H_0^1(\Omega) \times L^2(\omega)} = \sqrt{\int_\Omega |\nabla(I_M - P_M)\lambda^*|^2\, dx + \int_\omega |(I_O - P_O)f^*|^2\, dx}.$$

The two partial derivatives $\partial_\lambda G_P^\delta \in H_0^1(\Omega)$ and $\partial_f G_P^\delta \in L^2(\omega)$ of the Fréchet derivative of $G_P^\delta$ at point $(\lambda^*, f^*)$ are given by the weak formulation

$$\int_\Omega \nabla(\partial_\lambda G_P^\delta) \cdot \nabla\mu\, dx = \int_\Omega \nabla u^* \cdot \nabla\mu\, dx + \frac{\delta}{\|(I - P)(\lambda^*, f^*)\|_{H_0^1(\Omega) \times L^2(\omega)}} \int_\Omega \nabla(I_M - P_M)\lambda^* \cdot \nabla\mu\, dx, \quad \forall \mu \in H_0^1(\Omega) \quad (31)$$

and the formula

$$\partial_f G_P^\delta = u^*|_\omega + \frac{\delta}{\|(I - P)(\lambda^*, f^*)\|_{H_0^1(\Omega) \times L^2(\omega)}} (I_O - P_O)f^* - f^\delta. \quad (32)$$

## 4.2 The Cauchy problem for the Laplace equation

We again consider a bounded Lipschitz domain $\Omega \subset \mathbb{R}^d$, $d > 1$, such that its boundary $\partial\Omega$ is partitioned into two sets $\Gamma$ and $\tilde{\Gamma}$. More precisely, $\Gamma$ and $\tilde{\Gamma}$ are non empty open sets for the topology induced on $\partial\Omega$ from the topology on $\mathbb{R}^d$, that is $\partial\Omega = \overline{\Gamma} \cup \overline{\tilde{\Gamma}}$ and $\Gamma \cap \tilde{\Gamma} = \emptyset$.

The Cauchy problem consists here, for some data $f \in L^2(\Gamma)$, in finding $u \in H^1(\Omega)$ such that

$$\begin{cases} \Delta u &= 0 \quad \text{in } \Omega \\ u &= f \quad \text{on } \Gamma \\ \partial_\nu u &= 0 \quad \text{on } \Gamma, \end{cases} \tag{33}$$

where $\nu$ is the outward unit normal to $\partial\Omega$. This kind of problem arises when we have $\partial_\nu u = 0$ and measure $u = f$ on some accessible part $\Gamma$ of the boundary of the structure, while the complementary part $\tilde{\Gamma}$ of the boundary is not accessible. In practice those measurements are contaminated by some noise. Due to Holmgren's theorem, the Cauchy problem (33) has at most one solution. However it is ill-posed in the sense of Hadamard: existence may not hold for some data $f$, as for example shown in [1].

Let us show that the problem (33) is a particular case of the general framework exposed in section 3. It corresponds to spaces $V = H^1(\Omega)$, $M = \{\lambda \in H^1(\Omega), \lambda|_{\tilde{\Gamma}} = 0\}$, $O = L^2(\Gamma)$, the bilinear form $b$ on $V \times M$ such that $b(u, \lambda) = \int_\Omega \nabla u \cdot \nabla \lambda \, dx$, and the associated operator $B : V \to M$ such that $(Bu, \lambda)_M = b(u, \lambda)$ for all $(u, \lambda) \in V \times M$, while the operator $C : H^1(\Omega) \to L^2(\Gamma)$ is the trace operator. Due to Poincaré inequality, the space $M$ can be equipped with the $H^1$ semi-norm.

**Lemma 4.4.** *For $f \in L^2(\Gamma)$, $u \in H^1(\Omega)$ satisfies problem (33) if and only if $Bu = 0$ and $Cu = f$.*

*Proof.* The solution $u \in H^1(\Omega)$ satisfies problem (33) iff $u|_\Gamma = f$ and

$$b(u, \mu) = \int_\Omega \nabla u \cdot \nabla \mu \, dx = 0, \quad \forall \mu \in H^1(\Omega) \quad \text{such that} \quad \mu|_{\tilde{\Gamma}} = 0,$$

in other words iff $Bu = 0$ and $Cu = f$. $\square$

The injectivity of the operator $A = (B, C)$ is a direct consequence of uniqueness for the Cauchy problem. Hence we have checked that all properties of the general framework exposed in section 3 are satisfied. The corresponding formulation (16) is: for $f \in L^2(\Gamma)$, find $(u_\varepsilon, \lambda_\varepsilon) \in H^1(\Omega) \times \{\lambda \in H^1(\Omega), \lambda|_{\tilde{\Gamma}} = 0\}$ such that for all $(v, \mu) \in H^1(\Omega) \times \{\lambda \in H^1(\Omega), \lambda|_{\tilde{\Gamma}} = 0\}$,

$$\begin{cases} \varepsilon \int_\Omega (u_\varepsilon v + \nabla u_\varepsilon \cdot \nabla v) \, dx + \int_\Gamma u_\varepsilon v \, ds + \int_\Omega \nabla v \cdot \nabla \lambda_\varepsilon \, dx = \int_\Gamma f v \, ds, \\ \int_\Omega \nabla u_\varepsilon \cdot \nabla \mu \, dx - \int_\Omega \nabla \lambda_\varepsilon \cdot \nabla \mu \, dx = 0. \end{cases} \tag{34}$$

As in the previous application, by a Poincaré-type inequality the mapping $v \mapsto (\int_\Omega |\nabla v|^2 \, dx + \int_\Gamma v^2 \, ds)^{1/2}$ is a norm which is equivalent to the standard norm on $H^1(\Omega)$. Then thanks to Remark 3.2, the term $\int_\Omega u_\varepsilon v \, dx$ in (34) could be dropped. With a view to applying Corollary 3.1, let us identify $(\text{Range } A)^\perp$. Contrary to the previous data assimilation case, such space is reduced to $\{0\}$.

**Lemma 4.5.** *The operator $A$ has a dense range, or equivalently, $(\text{Range } A)^\perp = \{0\}$.*

*Proof.* Let us assume that $(\mu, h) \in M \times O = \{\lambda \in H^1(\Omega), \lambda|_{\tilde{\Gamma}} = 0\} \times L^2(\Gamma)$ satisfy, for all $v \in V = H^1(\Omega)$,

$$(Bv, \mu)_M + (Cv, h)_O = 0,$$

that is

$$\int_\Omega \nabla v \cdot \nabla \mu \, dx + \int_\Gamma v h \, ds = 0.$$

Choosing $v = \varphi \in C_0^\infty(\Omega)$ implies that $\Delta \mu = 0$ in the sense of distributions. Then, by the Green formula, we have that for all $v \in H^1(\Omega)$,

$$\langle v, \partial_\nu \mu \rangle_{H^{1/2}(\partial\Omega), H^{-1/2}(\partial\Omega)} + \int_\Gamma v h \, ds = 0,$$

where $\langle \cdot, \cdot \rangle_{H^{1/2}(\partial\Omega), H^{-1/2}(\partial\Omega)}$ denotes duality pairing between $H^{1/2}(\partial\Omega)$ and $H^{-1/2}(\partial\Omega)$. Considering $\tilde{h}$ as the extension of $h$ by 0 on $\partial\Omega$, we get that for all $v \in H^1(\Omega)$,

$$\langle v, \partial_\nu \mu + \tilde{h} \rangle_{H^{1/2}(\partial\Omega), H^{-1/2}(\partial\Omega)} = 0,$$

hence $\partial_\nu \mu + \tilde{h} = 0$ on $\partial\Omega$, that is $\partial_\nu \mu = 0$ on $\tilde{\Gamma}$ and $\partial_\nu \mu + h = 0$ on $\Gamma$. From uniqueness of the Cauchy problem applied to $\mu$, we conclude that $\mu = 0$ in $\Omega$, which implies in turn that $h = 0$ in $\Gamma$. This completes the proof. $\square$

We are now in a position to apply Corollary 3.1 in the simple case when $(\text{Range}\, A)^\perp = \{0\}$. Let us consider some data $f^\delta \in L^2(\Gamma)$ such that $\|f^\delta\|^2_{L^2(\Gamma)} > \delta$. There exists a unique $\varepsilon > 0$ such that

$$\int_\Omega |\nabla \lambda_\varepsilon^\delta|^2 \, dx + \int_\Gamma |u_\varepsilon^\delta - f^\delta|^2 \, ds = \delta^2,$$

where $(u_\varepsilon^\delta, \lambda_\varepsilon^\delta)$ is the unique solution to the problem (34) for data $f^\delta$.

**Remark 4.3.** *As a conclusion of sections 4.1 and 4.2, strictly speaking the Morozov's principle is easier to apply in the case of the Cauchy problem than in the case of the data assimilation problem for the Laplace equation. Indeed, the underlying operator $A$ has a dense range for the Cauchy problem while it has not for the data assimilation problem. As far as we know, such fact has never been highlighted so far. However, in view of Remark 2.2, the additional condition to check when $A$ has not a dense range is very unlikely to be violated.*

## 4.3 The data assimilation problem for the heat equation

Let us consider again some domains $\omega \Subset \Omega$ having the same properties as in section 4.1. Let us introduce $T > 0$, as well as $Q = \Omega \times (0, T)$, $q = \omega \times (0, T)$, $\Sigma = \partial\Omega \times (0, T)$ and $\sigma = \partial\omega \times (0, T)$. The data assimilation problem for the heat equation consists, for some data $f \in L^2(q)$, in finding $u \in L^2(0, T; H^1(\Omega))$ such that $\partial_t u - \Delta u = 0$ in $Q$ (in the sense of distributions) and $u|_q = f$. As for the Laplacian case, such problem has at most one solution in virtue of Holmgren's theorem but is ill-posed.
In view of the general framework exposed in section 3, the data assimilation problem for the heat equation corresponds to spaces $V = L^2(0, T; H^1(\Omega))$, $M = H_0^1(Q)$, $O = L^2(q)$, the bilinear form $b$ on $L^2(0, T; H^1(\Omega)) \times M$ such that $b(u, \lambda) = \int_Q (-u \, \partial_t \lambda + \nabla u \cdot \nabla \lambda) \, dxdt$, and the associated operator $B : L^2(0, T; H^1(\Omega)) \to H_0^1(Q)$ such that $(Bu, \lambda)_M = b(u, \lambda)$ for all $(u, \lambda) \in L^2(0, T; H^1(\Omega)) \times H_0^1(Q)$, while the operator $C : L^2(0, T; H^1(\Omega)) \to L^2(q)$ is the restriction operator. The space $H_0^1(Q)$ is equipped with the $H^1$ semi-norm in $Q$. We have the following lemma, the proof of which is very similar to the proof of Lemma 4.1.

**Lemma 4.6.** *For $f \in L^2(q)$, $u \in L^2(0, T; H^1(\Omega))$ satisfies the data assimilation problem for the heat equation if and only if $Bu = 0$ and $Cu = f$.*

The injectivity of the operator $A = (B, C)$ is a direct consequence of Holmgren's Theorem. The corresponding formulation (16) is: for $f \in L^2(q)$, find $(u_\varepsilon, \lambda_\varepsilon) \in L^2(0, T; H^1(\Omega)) \times H_0^1(Q)$ such that for all $(v, \mu) \in L^2(0, T; H^1(\Omega)) \times H_0^1(Q)$,

$$\begin{cases} \varepsilon \int_Q (u_\varepsilon \, v + \nabla u_\varepsilon \cdot \nabla v) \, dxdt + \int_q u_\varepsilon \, v \, dxdt + \int_Q (-v \, \partial_t \lambda_\varepsilon + \nabla v \cdot \nabla \lambda_\varepsilon) \, dxdt = \int_q f \, v \, dxdt, \\ \int_Q (-u_\varepsilon \, \partial_t \mu + \nabla u_\varepsilon \cdot \nabla \mu) \, dxdt - \int_Q (\partial_t \lambda_\varepsilon \, \partial_t \mu + \nabla \lambda_\varepsilon \cdot \nabla \mu) \, dxdt = 0. \end{cases} \quad (35)$$

It remains to identify $(\text{Range}\, A)^\perp$.

**Lemma 4.7.** *Denoting $H_q = H_0^1(0, T; L^2(\omega)) \cap L^2(0, T; H_0^2(\omega))$, we have that*

$$(\text{Range}\, A)^\perp = \{(\mu, h) \in H_q \times L^2(q), \quad \partial_t \mu + \Delta \mu = h\}.$$

*Proof.* Let us consider $(\mu, h) \in H_0^1(Q) \times L^2(\omega)$. We have $(\mu, h) \in (\text{Range}\, A)^\perp$ if and only if for all $v \in L^2(0, T; H^1(\Omega))$

$$(Bv, \mu)_M + (Cv, h)_O = 0,$$

that is
$$\int_Q (-v\,\partial_t \mu + \nabla v \cdot \nabla \mu)\,dxdt + \int_q v\,h\,dxdt = 0,$$
which is equivalent to
$$\begin{cases} \partial_t \mu + \Delta \mu &= 1_q\, h \quad \text{in } Q \\ \partial_\nu \mu &= 0 \quad \text{on } \Sigma. \end{cases}$$

We observe that $(\mu, h) \in (\text{Range } A)^\perp$ implies that $\mu$ satisfies $\partial_t \mu + \Delta \mu = 0$ in $Q \setminus \bar{q}$, $\mu = 0$ and $\partial_\nu \mu = 0$ on $\Sigma$, hence $\mu = 0$ in $Q \setminus \bar{q}$ from uniqueness of the lateral Cauchy problem for the heat equation. We obtain that $\mu = 0$ and $\partial_\nu \mu = 0$ on $\sigma$. As a result, for $(\mu, h) \in (\text{Range } A)^\perp$, we have $\mu \in H_0^1(0, T; L^2(\omega)) \cap L^2(0, T; H_0^2(\omega))$ and $\partial_t \mu + \Delta \mu = h$ in $q$. The converse statement is straightforward. □

Let us specify the orthogonal projection of $(0, f) \in H_0^1(Q) \times L^2(q)$ on $(\text{Range } A)^\perp$.

**Lemma 4.8.** *The projection of $(0, f) \in H_0^1(Q) \times L^2(\omega)$ on $(\text{Range } A)^\perp$ is the pair $(\lambda_\perp, f_\perp) \in H_0^1(Q) \times L^2(\omega)$, where $\lambda_\perp$ is the unique solution in $H_q$ of the weak formulation set in the subdomain $q$:*

$$\int_q (\partial_t \lambda_\perp + \Delta \lambda_\perp)(\partial_t \mu + \Delta \mu)\,dxdt + \int_q (\partial_t \lambda_\perp\,\partial_t \mu + \nabla \lambda_\perp \cdot \nabla \mu)\,dxdt = \int_q f\,(\partial_t \mu + \Delta \mu)\,dxdt, \quad \forall \mu \in H_q, \quad (36)$$

*and $f_\perp = \partial_t \lambda_\perp + \Delta \lambda_\perp$.*

The proof is omitted since it is similar to the proof of Lemma 4.3. Applying Corollary 3.1 to our problem, we get that for $f^\delta \in L^2(q)$ and $(\lambda_\perp^\delta, f_\perp^\delta)$ the orthogonal projection of $(0, f^\delta)$ on $(\text{Range } A)^\perp$, which can be computed with the help of Lemma 4.8, if we assume that

$$\|\lambda_\perp^\delta\|_{H_0^1(Q)}^2 + \|f_\perp^\delta\|_{L^2(q)}^2 < \delta^2 < \|f^\delta\|_{L^2(q)}^2,$$

then there exists a unique $\varepsilon > 0$ such that

$$\int_Q \left((\partial_t \lambda_\varepsilon)^2 + |\nabla \lambda_\varepsilon^\delta|^2\right)\,dxdt + \int_q |u_\varepsilon^\delta - f^\delta|^2\,dxdt = \delta^2,$$

where $(u_\varepsilon^\delta, \lambda_\varepsilon^\delta)$ is the unique solution to the problem (35) for data $f^\delta$.

## 5 Numerical experiments

All our numerical experiments concern the data assimilation problem for the Laplace equation which is addressed in section 4.1, in the two-dimensional case ($d = 2$). We have used the Freefem library [19] for all our finite element computations. The domain $\Omega$ is the square $(0, 1) \times (0, 1)$, while for $\omega$ we consider three different domains, all having the same surface $|\omega| = 0.4\,|\Omega|$, and which are represented on figure 1. The domain 1 is the delicate case of interior data, the domain 2 is the easier case of exterior data, while the domain 3 is an intermediate case where the data are distributed.

### 5.1 Illustration of Theorem 2.1

We first want to show a numerical illustration of Theorem 2.1 in the case of the data assimilation problem for the Laplace equation exposed in section 4.1. Let us consider the harmonic function $u(x, y) = 1 - x^3 + 3xy^2$ and $f = u|_\omega$, for $\omega$ corresponding to domain 2 in figure 1. For a given amplitude of noise $\delta$, the first thing to do is, in view of assumptions (1) and (7), to produce some noisy data $g^\delta = (\ell^\delta, f^\delta) \in H_0^1(\Omega) \times L^2(\omega) = H$ such that $\|g_\perp^\delta\|_H < \delta < \|g^\delta\|_H$ and $\|g^\delta - g\|_H = \delta$, where $g = (0, f)$. Here we choose $\delta = 0.1$. Note that in our data assimilation abstract framework, we have assumed that $\ell^\delta = 0$. However, with a view to building compatible noisy data more easily, for our validation of Theorem 2.1 we tolerate that $\ell^\delta \neq 0$. The consequence is that the weak formulation (24) becomes: find $(u_\varepsilon^\delta, \lambda_\varepsilon^\delta) \in H^1(\Omega) \times H_0^1(\Omega)$ such that for all $(v, \mu) \in H^1(\Omega) \times H_0^1(\Omega)$,

$$\begin{cases} \varepsilon \int_\Omega (u_\varepsilon^\delta v + \nabla u_\varepsilon^\delta \cdot \nabla v)\,dx + \int_\omega u_\varepsilon^\delta v\,dx + \int_\Omega \nabla v \cdot \nabla \lambda_\varepsilon^\delta\,dx = \int_\omega f^\delta v\,dx, \\ \int_\Omega \nabla u_\varepsilon^\delta \cdot \nabla \mu\,dx - \int_\Omega \nabla \lambda_\varepsilon^\delta \cdot \nabla \mu\,dx = \int_\Omega \nabla \ell^\delta \cdot \nabla \mu\,dx. \end{cases} \quad (37)$$

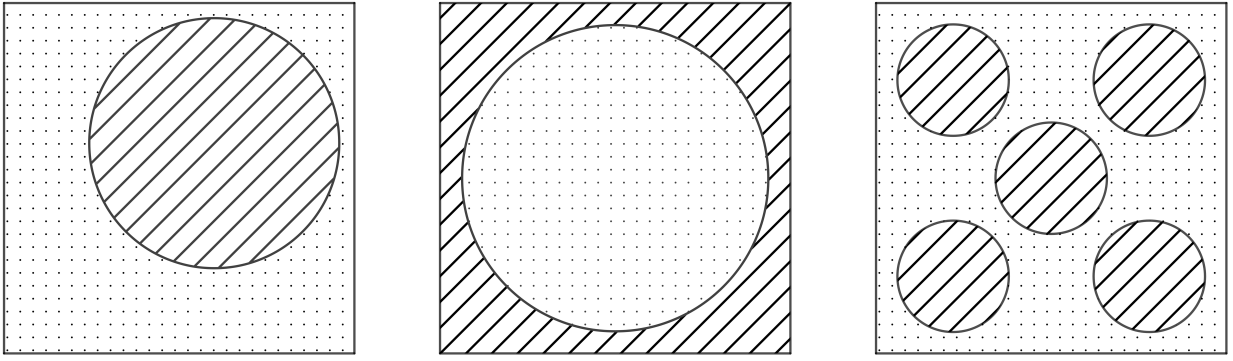

Figure 1: Three different domains $\omega$ (the domain $\omega$ is hatched). Left (domain 1): a disk. Middle (domain 2): the exterior of a disk. Right (domain 3): 5 small disks.

Note that contrary to (24), the right-hand side of the second equation in (37) is not 0 any more. On the one hand we introduce the function $w = (x^2+y^2)/4$ and set $(\ell_{/\!/}, f_{/\!/}) = (\Delta w, w|_\omega)$. On the other hand we introduce the indicator function $\tilde{f}(x,y) = 1$ if $x^2 + y^2 > 1$ and $\tilde{f}(x,y) = 0$ otherwise. For such $\tilde{f} \in L^2(\omega)$, let us introduce $(\ell_\perp, f_\perp) \in H_0^1(\Omega) \times L^2(\omega)$ the orthogonal projection of $(0, \tilde{f})$ on $(\text{Range}\, A)^\perp$ given by Lemma 4.3 and more precisely by the weak formulation (26). We have hence obtained a pair $g_{/\!/} = (\ell_{/\!/}, f_{/\!/}) \in \text{Range}\, A$ and a pair $g_\perp = (\ell_\perp, f_\perp) \in (\text{Range}\, A)^\perp$. The idea is to search some noisy data $g^\delta = (\ell^\delta, f^\delta)$ in the form

$$g^\delta = g + \alpha g_{/\!/} + \beta g_\perp,$$

where $\alpha, \beta > 0$ are uniquely defined such that

$$\|g_\perp^\delta\|_H = \frac{\delta}{2}, \quad \|g^\delta - g\|_H = \delta.$$

A straightforward computation yields

$$\|g^\delta\|_H^2 = \|g\|_H^2 + \delta^2 + 2\alpha(f, f_{/\!/})_{L^2(\omega)}.$$

Observing that $(f, f_{/\!/})_{L^2(\omega)} \geq 0$, we conclude that $\|g^\delta\|_H > \delta$, that is our artificial noisy data $g^\delta$ satisfies both assumptions (1) and (7). Using this noisy data $g^\delta$, for $\varepsilon > 0$ we compute

$$E^\delta(\varepsilon) = \|Au_\varepsilon^\delta - g^\delta\|_H = \sqrt{\int_\Omega |\nabla \lambda_\varepsilon^\delta|^2 \, dx + \int_\omega |u_\varepsilon^\delta - f^\delta|^2 \, dx},$$

where $(u_\varepsilon^\delta, \lambda_\varepsilon^\delta) \in H^1(\Omega) \times H_0^1(\Omega)$ is the solution to problem (37). The computation of $(u_\varepsilon^\delta, \lambda_\varepsilon^\delta)$ is based on a Finite Element Method, the space $H^1(\Omega)$ being approximated with the help of the classical $P1$ finite elements based on triangles. The resolution of problem (26) is also based on a FEM, the space $H^2(\omega)$ being approximated by the $C^1$ conforming Hsieh-Clough-Tocher triangular finite elements (see for example [13] for a description of such element). In both cases, the mesh size is $h = 1/20$. The graph of the function $E^\delta$ of $\varepsilon$ is plotted on figure 2. We observe that the function $E^\delta$ is non increasing with $E^\delta(0^+) = \|g_\perp^\delta\|_H = \delta/2$ and $E^\delta(+\infty) = \|g^\delta\|_H > \delta$, in accordance with Theorem 2.1, which seems to indicate that the chosen mesh size $h$ is sufficiently small. The important conclusion is that the application of the Morozov's principle seems to be relevant even for the discretized problem, and not only for the continuous one.

### 5.2 Validation of the duality method

After having observed that the computation of the Morozov's solution associated with the regularized problem (24) actually makes sense after discretization, it is now natural to test the duality method to obtain such solution, which consists in minimizing the functional $G_0^\delta$ given by (30) and then applying the operator $A^*$ to the obtained solution. Again, the Finite Element Method is used, the space $H^1(\Omega)$ being approximated by $P1$ finite elements, while the space $L^2(\omega)$ is approximated with $P0$ finite elements. The minimization of $G_0^\delta$ is an iterative method based on the computations of the partial derivatives $\partial_\lambda G_0^\delta$ and $\partial_f G_0^\delta$ given by (31) and (32), respectively. Note that an alternative technique consists, in order to avoid minimizing the non-quadratic functional $G_0^\delta$, to find, for $n \geq 1$, the minimum $(\lambda_n^*, f_n^*) \in H_0^1(\Omega) \times L^2(\omega)$ of the quadratic functional

$$H_n^\delta(\lambda^*, f^*) = \frac{1}{2} \|u^*(\lambda^*, f^*)\|_{H^1(\Omega)}^2 + \varepsilon_n \|(\lambda^*, f^*)\|_{H_0^1(\Omega) \times L^2(\omega)}^2 - \int_\omega f^\delta f^* \, dx, \qquad (38)$$

where $u^*$ is the solution to problem (28) and where $\varepsilon_n = \delta / \|(\lambda_{n-1}^*, f_{n-1}^*)\|_{H_0^1(\Omega) \times L^2(\omega)}$. This technique was introduced in [16], where it is proved that in the case when the operator $A$ has a dense range, the sequence formed by the minima of the quadratic problems converges to the unique minimum of the non-quadratic one. Despite we have not extended this result to the case when $A$ does not have a dense range, which is the case with our data assimilation problem for the Laplace equation, such alternative technique seems to work in the sense that, at the numerically level, the sequence $(\lambda_n^*, f_n^*)_n$ converges to the same minimizer as the one obtained by directly minimizing the functional $G_0^\delta$ in $H_0^1(\Omega) \times L^2(\omega)$ given by (30).

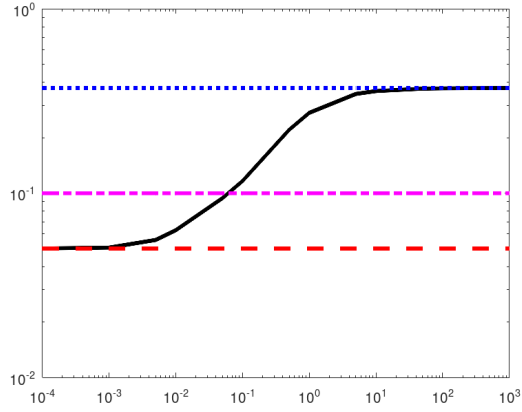

Figure 2: Graph of the function $E^\delta$ of $\varepsilon$. The horizontal lines correspond to the values $\|g^\delta_\perp\|_H = \delta/2$, $\|g^\delta - g\|_H = \delta$ and $\|g^\delta\|_H > \delta$, respectively.

In what follows, we will consider two kinds of harmonic function $u$ given either by

$$u(x,y) = \alpha_\infty (1 - x^3 + 3xy^2) \qquad (39)$$

or

$$u(x,y) = \alpha_\infty \sin(4(1-y))e^{4x}, \qquad (40)$$

where $\alpha_\infty$ is calibrated such that $\|u\|_{L^\infty(\Omega)} = 1$ and set $f = u|_\omega$, for $\omega$ one of the three domains described in figure 1. We design a noisy data $f^\delta$ by adding a pointwise random quantity to $f$ in such a way that

$$\|f^\delta - f\|_{L^2(\omega)} = \delta_r \|f\|_{L^2(\omega)} = \delta,$$

where $\delta_r$ is some prescribed relative noise, that is $\delta_r = 2\%$, $\delta_r = 5\%$ or $\delta_r = 10\%$.

In figure 3, we have plotted the exact solution $u$ given (39). In the figures 4, 5 and 6 below and for $\delta_r = 10\%$, we have plotted the solution $u^\delta$ obtained from the minimization of $G_0^\delta$ (the so-called Morozov's solution) and their difference $u - u^\delta$, for $\omega$ being either of the domains 1, 2 and 3, respectively.

In figure 7, we have plotted the exact solution $u$ given (40). In figures 8, 9 and 10 below and for $\delta_r = 10\%$, we have again plotted the Morozov's solution $u^\delta$ and the difference $u - u^\delta$, for the three different domains $\omega$.

From these numerical results, for both exact solutions given by (39) and (40), we observe that the reconstruction is the worst for domain 1, the best for domain 2 and intermediate for domain 3.

In figure 11, for the exact solution (40), the domain 3 and $\delta_r = 10\%$, we have plotted the error $\|u_\varepsilon^\delta - u\|_{H^1(\Omega)}$ as a function of $\varepsilon$, where $(u_\varepsilon^\delta, \lambda_\varepsilon^\delta)$ is the solution to problem (24) for data $f^\delta$. On the same graph, the value of $\varepsilon(\delta)$ corresponding to the Morozov choice, which is associated to the Morozov's solution $u^\delta$, is represented.

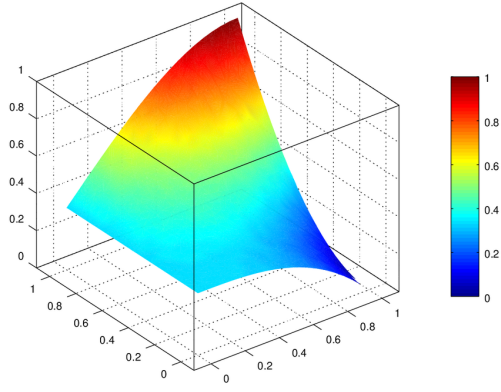

Figure 3: Exact solution (39).

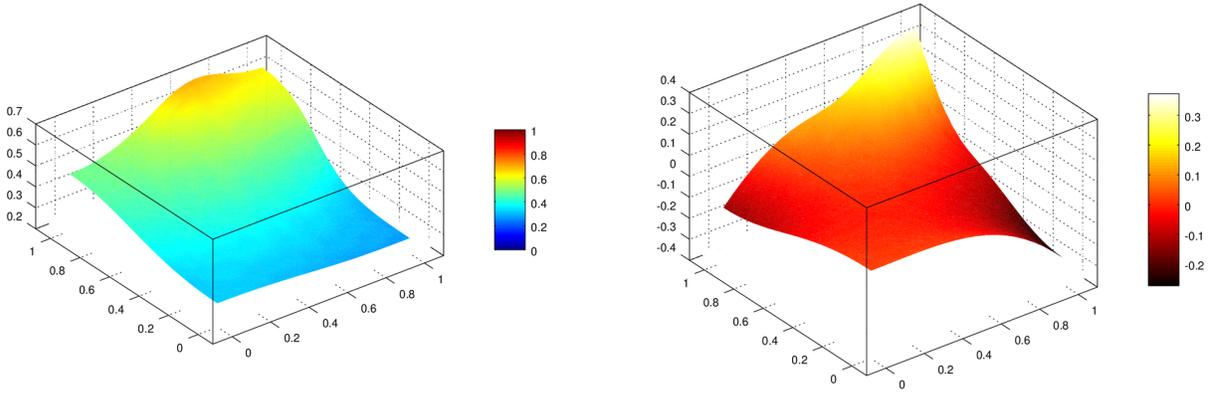

Figure 4: Comparison of the Morozov's solution and the exact solution (39) for domain 1. Left: function $u^\delta$. Right: function $u - u^\delta$. Relative error is $\|u - u^\delta\|_{H^1(\Omega)}/\|u\|_{H^1(\Omega)} = 0.55$

We clearly see that for $\varepsilon < \varepsilon(\delta)$, the error reaches a plateau, which indicates that the Morozov choice for $\varepsilon$ is a relevant value. Indeed, from a numerical point of view, selecting a value of $\varepsilon$ smaller than $\varepsilon(\delta)$ would deteriorate the condition number of the matrix to invert without improving the quality of the reconstruction in terms of the error $\|u^\delta_\varepsilon - u\|_{H^1(\Omega)}$. In addition, the presence of a plateau instead of a clear minimum is a consequence of the discretization. In practice we observe that the function $\delta \mapsto \varepsilon(\delta)$ is non decreasing (see for example [6] for a similar case). Such monotonicity property cannot be proved theoretically without making assumptions on the mapping $\delta \mapsto f^\delta \in L^2(\omega)$.

Now let us analyze the influence of the surface of the domain $\omega$, that is $|\omega| = \alpha|\Omega|$, with $\alpha = 1/5$, $\alpha = 1/4$, $\alpha = 1/3$ and $\alpha = 2/5$, for $\omega$ being the domain 3, $u$ given by (39) and $\delta_r = 5\%$. The obtained results are given in the table 12 below. As expected, the larger is the measure of $\omega$, the better is the quality of the identification.

We complete this numerical section by testing the minimization of the functional (30) when the projector $P = (P_M, P_O)$ is different from 0. Let us first consider an example of projector $P_O$, that is the operator $P_O : L^2(\omega) \to L^2(\omega)$ such that

$$P_O f = \frac{1}{|\omega|} \left( \int_\omega f \, dx \right) 1_\omega, \quad \forall f \in L^2(\omega).$$

This choice is motivated by the following remark: let us assume that for a given $x \in \omega$, the value of $(f - f^\delta)(x)$ is a random function of 0 mean value. By ergodicity, computing the mean value of $f^\delta(x)$ at point $x$ for all realization amounts to compute, for a given realization, the mean value of $f^\delta$ on $\omega$. This is why we impose

$$\int_\omega u^\delta \, dx = \int_\omega f^\delta \, dx$$

for the regularized solution $u^\delta$, that is $P_O C u^\delta = P_O f^\delta$. Let us secondly consider an example of operator $P_M : H_0^1(\Omega) \to H_0^1(\Omega)$. We remark that the regularized solution $u^\delta$ to the problem (24) does not satisfy the equation $\Delta u^\delta = 0$ exactly, contrary to the exact solution $u$. Let us introduce the eigenvalues $\lambda_n$ and the corresponding normalized eigenfunctions $\varphi_n$ ($\|\varphi_n\|_{L^2(\Omega)} = 1$) of the Dirichlet-Laplacian operator $-\Delta$ in $\Omega$ for $n \in \mathbb{N}$, that is $\varphi_n \in H_0^1(\Omega)$ and $-\Delta\varphi_n = \lambda_n \varphi_n$ in $\Omega$. We also denote, for some $N \in \mathbb{N}$, the finite set $I_N = \{0, \cdots, N\}$ and the subspace of $M$ defined by $M_N = \text{Span}\{\varphi_n, n \in I_N\}$. We assume that $P_M$ is the orthogonal projector on $M_N$. On the one hand, we observe that $Bu = 0$ is equivalent to

$$\int_\Omega \nabla u \cdot \nabla \mu \, dx = 0, \ \forall \mu \in H_0^1(\Omega) \iff \int_\Omega \nabla u \cdot \nabla \varphi_n \, dx = 0, \ \forall n \in \mathbb{N} \iff \Delta u = 0.$$

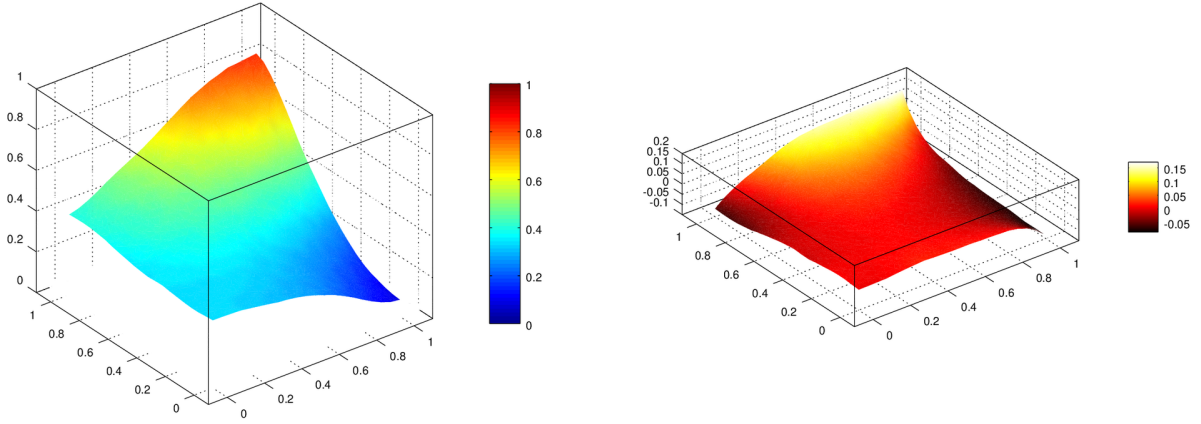

Figure 5: Comparison of the Morozov's solution and the exact solution (39) for domain 2. Left: function $u^\delta$. Right: function $u - u^\delta$. Relative error is $\|u - u^\delta\|_{H^1(\Omega)}/\|u\|_{H^1(\Omega)} = 0.33$

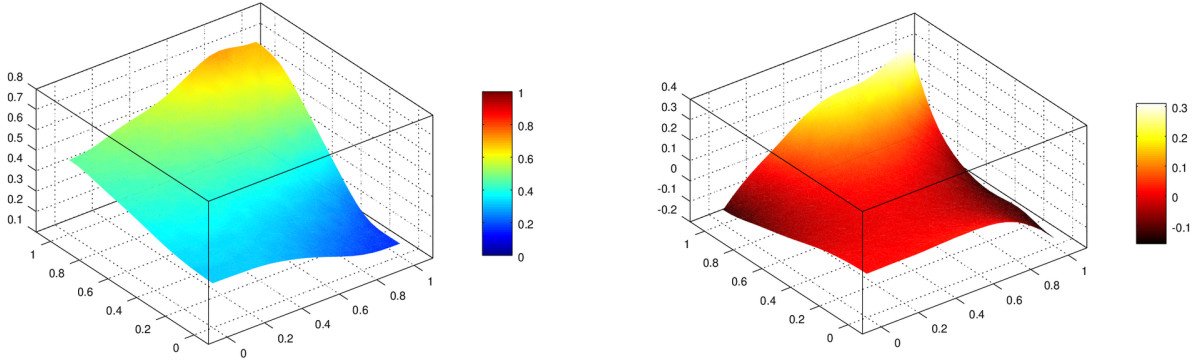

Figure 6: Comparison of the Morozov's solution and the exact solution (39) for domain 3. Left: function $u^\delta$. Right: function $u - u^\delta$. Relative error is $\|u - u^\delta\|_{H^1(\Omega)}/\|u\|_{H^1(\Omega)} = 0.52$

On the other hand, imposing $P_M B u^\delta = 0$ to the regularized solution $u^\delta$ exactly means that

$$(Bu^\delta, \varphi_n)_{H_0^1(\Omega)} = 0, \ \forall n \in I_N \iff \int_\Omega \nabla u^\delta \cdot \nabla \varphi_n \, dx = 0, \ \forall n \in I_N \iff \langle \Delta u^\delta, \varphi_n \rangle_{H^{-1}(\Omega), H_0^1(\Omega)} = 0, \ \forall n \in I_N.$$

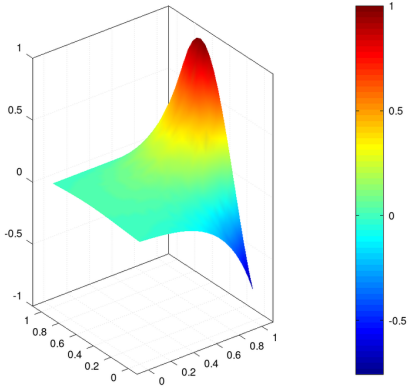

Figure 7: Exact solution (40).

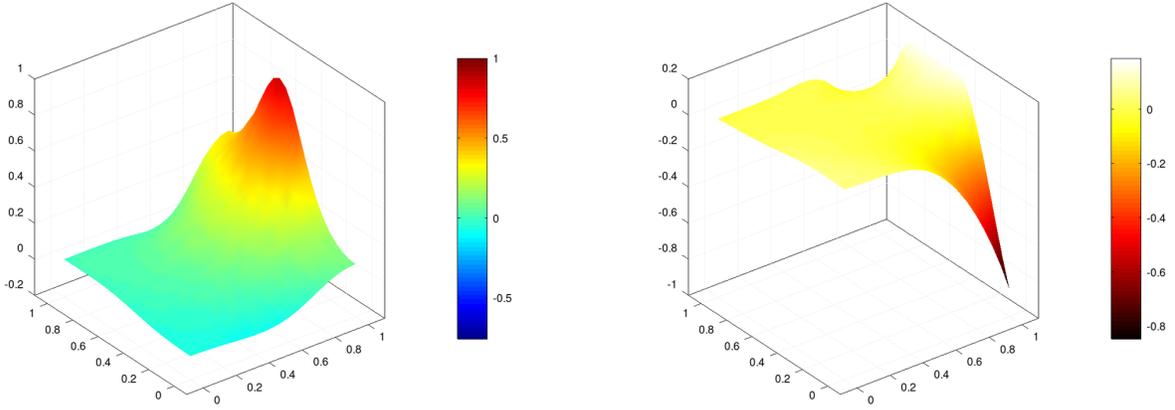

Figure 8: Comparison of the Morozov's solution and the exact solution (40) for domain 1. Left: function $u^\delta$. Right: function $u - u^\delta$. Relative error is $\|u - u^\delta\|_{H^1(\Omega)}/\|u\|_{H^1(\Omega)} = 0.60$

We conclude that $P_M B u^\delta = 0$ is equivalent to the fact that the restriction of the linear form $\Delta u^\delta$ to the subspace $M_N$ is null. In a sense, we have partially enforced the constraint $\Delta u^\delta = 0$. Based on the previous definitions of the operators $P_M$ and $P_O$, the operator $P = (P_M, P_O) : H_0^1(\Omega) \times L^2(\omega) \to H_0^1(\Omega) \times L^2(\omega)$ is an orthogonal projector which is compact. Unfortunately, such operator $P$ does not satisfy the required property Range $P \subset \overline{\text{Range } A}$. Indeed, proving that Range $P \subset \overline{\text{Range } A}$ amounts to prove that $(\lambda_\perp, f_\perp) = 0$ for any pair $(\lambda, f) \in \text{Range } P$, where $(\lambda_\perp, f_\perp) \in H_0^1(\Omega) \times L^2(\omega)$ is the orthogonal projection of $(\lambda, f)$ on $(\text{Range } A)^\perp$. The pair $(\lambda_\perp, f_\perp)$ is characterized by (27). Since $f$ is a constant function in $\omega$ and $\lambda = \sum_{n=0}^N (\lambda, \varphi_n)_{L^2(\Omega)} \varphi_n$ in $\Omega$, this characterization amounts to

$$\begin{cases} \lambda_\perp \in H_0^2(\omega) \quad \text{and} \quad f_\perp = \Delta \lambda_\perp, \\ \int_\omega \Delta \lambda_\perp \, \Delta \mu \, dx + \int_\omega \nabla \lambda_\perp \cdot \nabla \mu \, dx = \sum_{n=0}^N \lambda_n \left( \int_\Omega \lambda \, \varphi_n \, dx \right) \left( \int_\omega \varphi_n \, \mu \, dx \right), \quad \forall \mu \in H_0^2(\omega). \end{cases} \quad (41)$$

If $\lambda_\perp = 0$, we obtain from (41) that $\sum_0^N \mu_n \varphi_n = 0$ in $\omega$, denoting $\mu_n = \lambda_n (\lambda, \varphi)_{L^2(\Omega)}$. From the Lebeau-Robbiano spectral inequality shown in [23], we obtain that $\mu_n = 0$ for $n = 0, \cdots, N$, and we conclude that $\lambda = 0$. Hence $\lambda_\perp$ does not vanish unless $\lambda = 0$.

In order to restore the property Range $P \subset \overline{\text{Range } A}$, let us rather consider the eigenvalues $\lambda_n^0$ and the corresponding normalized eigenfunctions $\varphi_n^0$ of the Dirichlet-Laplacian operator $-\Delta$ in $\Omega \setminus \overline{\omega}$ for $n \in \mathbb{N}$, that is $\varphi_n^0 \in H_0^1(\Omega \setminus \overline{\omega})$

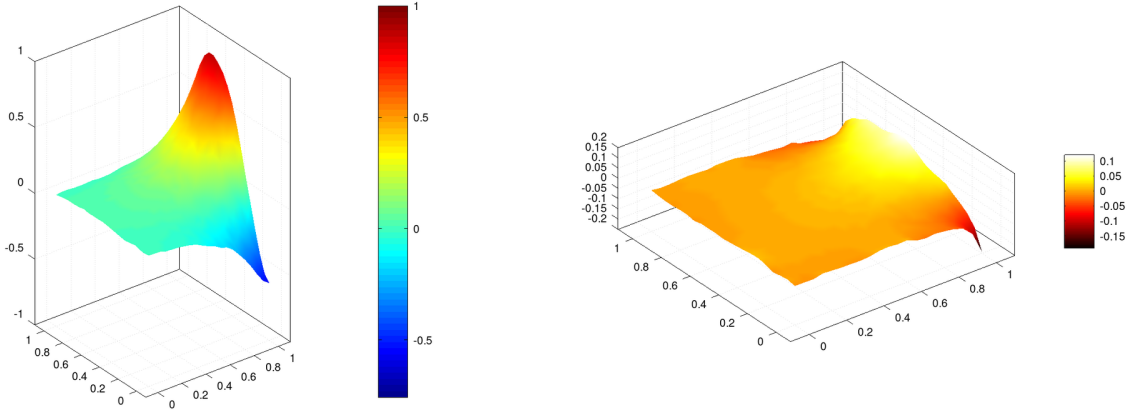

Figure 9: Comparison of the Morozov's solution and the exact solution (40) for domain 2. Left: function $u^\delta$. Right: function $u - u^\delta$. Relative error is $\|u - u^\delta\|_{H^1(\Omega)}/\|u\|_{H^1(\Omega)} = 0.20$

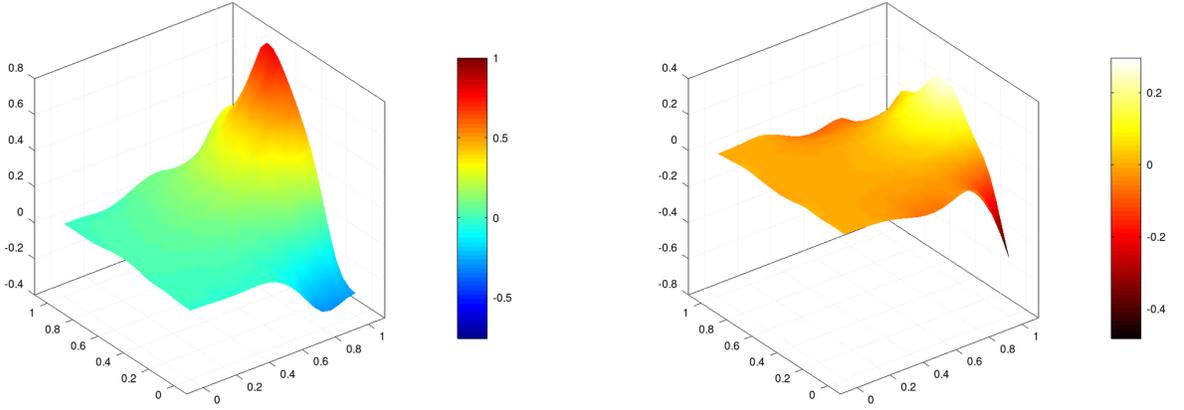

Figure 10: Comparison of the Morozov's solution and the exact solution (40) for domain 3. Left: function $u^\delta$. Right: function $u - u^\delta$. Relative error is $\|u - u^\delta\|_{H^1(\Omega)}/\|u\|_{H^1(\Omega)} = 0.40$

and $-\Delta \varphi_n^0 = \lambda_n^0 \varphi_n^0$ in $\Omega \setminus \overline{\omega}$. These functions $\varphi_n^0$, extended by 0 in $\omega$, are functions in $H_0^1(\Omega)$ which are still denoted $\varphi_n^0$, with a slight abuse of notation. The operator $P_M^0 : H_0^1(\Omega) \to H_0^1(\Omega)$ is then defined as the orthogonal projection of a function in $H_0^1(\Omega)$ on the $(N+1)$ first functions $\varphi_n^0$, which amounts to partially enforce the constraint $\Delta u^\delta = 0$ in $\Omega \setminus \overline{\omega}$ instead of $\Omega$. Finally, we form the operator $P^0 = (P_M^0, P_O) : H_0^1(\Omega) \times L^2(\omega) \to H_0^1(\Omega) \times L^2(\omega)$, which is a compact orthogonal projector. Since the functions $\varphi_n^0$ vanish in $\omega$, from (41) for $(\lambda_n, \varphi_n)$ replaced by $(\lambda_n^0, \varphi_n^0)$, we observe that $\text{Range}\, P^0 \subset \overline{\text{Range}\, A}$, which implies that $P^0$ is an admissible projector.

As a numerical experiment, we now minimize the functional (30) for the admissible projector $P^0 = (P_M^0, P_O)$ given above with $N = 4$, that is $Bu^\delta$ has null contributions on the first 5 eigenfunctions $\varphi_n^0$. The exact solution is given by (39), the domain $\omega$ corresponds to the domain 3 in the figure 1 with $|\omega| = 0.4 |\Omega|$, and we test three different values of the relative noise $\delta_r$, that is 2%, 5% and 10%. In the table 13 below, we compare how the constraints $P_O C u^\delta = P_O f^\delta$ and $P_M^0 B u^\delta = 0$ are satisfied with or without the projector $P^0$. More precisely, we compute $\int_\omega u^\delta \, dx$, which has to be compared to $\int_\omega f^\delta \, dx$, as well as $\|P_M^0 B u^\delta\|_{H_0^1(\Omega)}$.

The table 13 is complemented by the table 14 below, which provides the error between the exact and the Morozov's solution in the presence of the projector $P^0$. The line which corresponds to $\delta_r = 5\%$ in table 14 has to be compared with the line which corresponds to $\alpha = 2/5$ in table 12. There are almost the same. As a conclusion, using the projector $P^0$ in the functional (30) given above is an efficient way to impose some constraints on the solution $u^\delta$, in particular for large amplitudes of noise, without increasing the error between the regularized solution and the exact one.

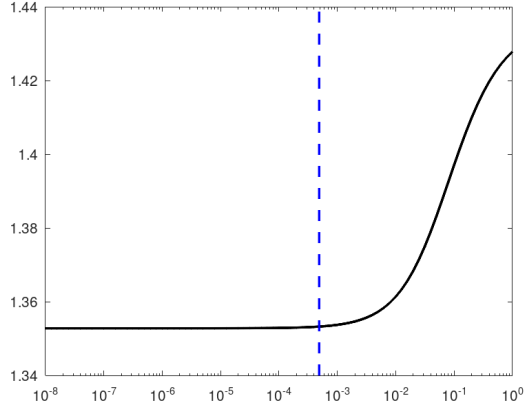

Figure 11: Error $\|u_\varepsilon^\delta - u\|_{H^1(\Omega)}$ as a function of $\varepsilon$. The blue line represents the Morozov value.

| $\alpha$ | $\|u - u^\delta\|_{H^1(\Omega)}$ | $\dfrac{\|u - u^\delta\|_{H^1(\Omega)}}{\|u\|_{H^1(\Omega)}}$ | $\|u - u^\delta\|_{L^2(\omega)}$ | $\dfrac{\|u - u^\delta\|_{L^2(\omega)}}{\|u\|_{L^2(\omega)}}$ |
|---|---|---|---|---|
| 1/5 | 0.50 | 0.55 | 0.02 | 0.10 |
| 1/4 | 0.47 | 0.52 | 0.02 | 0.092 |
| 1/3 | 0.41 | 0.46 | 0.018 | 0.072 |
| 2/5 | 0.3555 | 0.3909 | $1.6991\,10^{-2}$ | $6.0638\,10^{-2}$ |

Figure 12: Influence of $\alpha = |\omega|/|\Omega|$ on the error between the exact and the Morozov's solution

| | | $\int_\omega u^\delta\,dx$ | | $\|P_M^0 B u^\delta\|_{H_0^1(\Omega)}$ | |
|---|---|---|---|---|---|
| $\delta_r$ | $\int_\omega f^\delta\,dx$ | without projection | with projector $P^0$ | without projection | with projector $P^0$ |
| 2% | 0.1654 | 0.1649 | 0.1654 | $1.4016\,10^{-5}$ | $4.2373\,10^{-6}$ |
| 5% | 0.1653 | 0.1641 | 0.1653 | $1.2481\,10^{-5}$ | $2.5544\,10^{-6}$ |
| 10% | 0.1657 | 0.1617 | 0.1657 | $2.8569\,10^{-4}$ | $1.2073\,10^{-6}$ |

Figure 13: Comparison of how the constraints $P_O C u^\delta = P_O f^\delta$ and $P_M^0 B u^\delta = 0$ are satisfied, with and without projection

| $\delta_r$ | $\|u - u^\delta\|_{H^1(\Omega)}$ | $\dfrac{\|u - u^\delta\|_{H^1(\Omega)}}{\|u\|_{H^1(\Omega)}}$ | $\|u - u^\delta\|_{L^2(\omega)}$ | $\dfrac{\|u - u^\delta\|_{L^2(\omega)}}{\|u\|_{L^2(\omega)}}$ |
|---|---|---|---|---|
| 2% | 0.2630 | 0.2892 | $8.9105\,10^{-3}$ | $3.1800\,10^{-2}$ |
| 5% | 0.3573 | 0.3929 | $1.7038\,10^{-2}$ | $6.0805\,10^{-2}$ |
| 10% | 0.4774 | 0.5250 | $3.3145\,10^{-2}$ | $1.1829\,10^{-1}$ |

Figure 14: Influence of the projector $P^0$ on the error between the exact and the Morozov's solution

Lastly, it is tempting to produce the same tables 13 and 14 as before by using the more natural but non admissible projector $P = (P_M, P_O)$ instead of $P^0 = (P_M^0, P_O)$, that is tables 15 and 16 below. In view of the numerical results, and despite the operator $P$ is not admissible, it seems to be satisfactory from the point of view of its objectives: its enables us to enforce the constraints without deteriorating the error between the Morozov solution and the exact solution.

| $\delta_r$ | $\int_\omega f^\delta\,dx$ | $\int_\omega u^\delta\,dx$ | | $\|P_M B u^\delta\|_{H_0^1(\Omega)}$ | |
|---|---|---|---|---|---|
| | | without projection | with projector $P$ | without projection | with projector $P$ |
| 2% | 0.1654 | 0.1649 | 0.1654 | $3.2756\,10^{-4}$ | $1.7879\,10^{-4}$ |
| 5% | 0.1653 | 0.1641 | 0.1653 | $8.1457\,10^{-4}$ | $5.2360\,10^{-5}$ |
| 10% | 0.1657 | 0.1617 | 0.1657 | $2.5715\,10^{-3}$ | $1.8629\,10^{-4}$ |

Figure 15: Comparison of how the constraints $P_O C u^\delta = P_O f^\delta$ and $P_M B u^\delta = 0$ are satisfied, with and without projection

| $\delta_r$ | $\|u-u^\delta\|_{H^1(\Omega)}$ | $\dfrac{\|u-u^\delta\|_{H^1(\Omega)}}{\|u\|_{H^1(\Omega)}}$ | $\|u-u^\delta\|_{L^2(\omega)}$ | $\dfrac{\|u-u^\delta\|_{L^2(\omega)}}{\|u\|_{L^2(\omega)}}$ |
|---|---|---|---|---|
| 2% | 0.2648 | 0.2911 | $9.0449\,10^{-3}$ | $3.2280\,10^{-2}$ |
| 5% | 0.3573 | 0.3929 | $1.7056\,10^{-2}$ | $6.0871\,10^{-2}$ |
| 10% | 0.4767 | 0.5252 | $3.3119\,10^{-2}$ | $1.1820\,10^{-1}$ |

Figure 16: Influence of the projector $P$ on the error between the exact and the Morozov's solution

# Appendix: on duality in optimization

In this appendix, we give a brief review of the theory of duality in optimization exposed in [17]. We consider an optimization problem $(\mathscr{P})$ denoted the primal problem :

$$(\mathscr{P}) \quad \inf_{v\in V} F(v),$$

where $V$ is a Hilbert space, $F : V \to \overline{\mathbb{R}}$ a function $\neq +\infty$. Here we have denoted $\overline{\mathbb{R}} = \mathbb{R} \cup \{\pm\infty\}$. We recall the definition of the conjugate function, where $V^*$ denotes the dual space of $V$ and $\langle\cdot,\cdot\rangle_{V,V^*}$ denotes the duality bracket between $V$ and $V^*$.

**Definition 5.1.** *The congugate function $F^* : V^* \to \overline{\mathbb{R}}$ of $F$ is defined, for $u^* \in V^*$, by*

$$F^*(u^*) = \sup_{u\in V} \left(\langle u, u^*\rangle_{V,V^*} - F(u)\right).$$

Then we introduce the notion of perturbed problem. We consider a function $\Phi : V \times H \to \overline{\mathbb{R}}$, where $H$ is another Hilbert space, and $\Phi$ satisfies
$$\Phi(v,0) = F(v).$$
For all $q \in H$, we consider the perturbed problem $(\mathscr{P}_q)$:

$$(\mathscr{P}_q) \quad \inf_{v\in V} \Phi(v,q).$$

Then we define the dual problem of problem $(\mathscr{P})$ with respect to the perturbation $\Phi$. Let $\Phi^* : V^* \times H^* \to \overline{\mathbb{R}}$ be the conjugate function of $\Phi$. The dual problem, denoted $(\mathscr{P}^*)$, is the following optimization problem:

$$(\mathscr{P}^*) \quad \sup_{q^*\in H^*} -\Phi^*(0,q^*).$$

We have the following proposition.

**Proposition 5.1.**
$$(-\infty \leq) \quad \sup(\mathscr{P}^*) \leq \inf(\mathscr{P}) \quad (\leq +\infty).$$

*Proof.* For $q^* \in H^*$, we have
$$\Phi^*(0,q^*) = \sup_{u\in V,\,q\in H} \left(\langle u,0\rangle_{V,V^*} + \langle q,q^*\rangle_{H,H^*} - \Phi(u,q)\right),$$
so that, for all $u \in V$,
$$\Phi^*(0,q^*) \geq \langle 0,q^*\rangle_{H,H^*} - \Phi(u,0) = -\Phi(u,0).$$
We hence have $\forall u \in V$, $\forall q^* \in H^*$, $-\Phi^*(0,q^*) \leq \Phi(u,0)$, and finally $\sup(\mathscr{P}^*) \leq \inf(\mathscr{P})$. □

**Remark 5.1.** *Equality is not satisfied in general: when $\sup(\mathscr{P}^*) \neq \inf(\mathscr{P})$, we say that there is a duality gap.*

The following theorem, which is proved in [17], guarantees equality $\sup(\mathscr{P}^*) = \inf(\mathscr{P})$.

**Theorem 5.1.** *We assume that $\Phi$ is convex and that $\inf(\mathscr{P}) < +\infty$. If there exists $u_0 \in V$ such that $q \to \Phi(u_0, q)$ is finite and continuous at point $0$, then $\inf(\mathscr{P}) = \sup(\mathscr{P}^*) < +\infty$ and the problem $(\mathscr{P}^*)$ has solutions.*

Let us apply the above general theory to the primal problem:

$$\inf_{v \in K_0^\delta} L(v) \tag{42}$$

with

$$L(v) = \frac{1}{2}\|v\|_V^2, \quad K_0^\delta = \{v \in V \,;\, \|Av - g^\delta\|_H \leq \delta\},$$

where $A : V \to H$ is a bounded injective operator and $g^\delta \in H$. The problem (42) is indeed in the form $(\mathscr{P})$, that is

$$\inf_{v \in V} F(v),$$

if we define

$$F(u) = L(u) + \chi_{B_\delta}(Au),$$

and

$$\Phi(u, q) = L(u) + \chi_{B_\delta}(Au - q), \tag{43}$$

where $B_\delta \subset H$ is the closed ball of center $g^\delta$ and radius $\delta$, $\chi_{B_\delta}$ is the indicator function defined by

$$\begin{cases} \chi_{B_\delta}(q) = 0 & \text{if } q \in B_\delta \\ \chi_{B_\delta}(q) = +\infty & \text{if } q \notin B_\delta. \end{cases}$$

We have the following theorem, which is a consequence of Theorem 2.3 when $P = 0$.

**Theorem 5.2.** *For $g^\delta \in H$ satisfying assumption (1), the problem $(\mathscr{P})$ given by (42) has a unique solution $u^\delta \in V$.*

Now let us derive the dual problem $(\mathscr{P}^*)$.

**Proposition 5.2.** *The dual problem which corresponds to the primal problem (42) is equivalent to*

$$(\mathscr{P}^*) \quad \inf_{q^* \in H^*} G_0^\delta(q^*) = \inf_{q^* \in H} \left( \frac{1}{2}\|A^*q^*\|_{V^*}^2 + \delta \|q^*\|_{H^*} - \langle g^\delta, q^*\rangle_{H, H^*} \right). \tag{44}$$

*Proof.* Let us form the dual problem $(\mathscr{P}^*)$ which, after some simple computations, is defined by

$$(\mathscr{P}^*) \quad \sup_{q*\in H} -\Phi^*(0, q^*) = \sup_{q*\in H} \left( -L^*(A^*q^*) - \chi^*_{B_\delta}(-q^*) \right).$$

It remains to compute $L^*$ and $\chi^*_{B_\delta}$.
It is easy to see that

$$L^*(v^*) = \frac{1}{2}\|v^*\|_{V^*}^2.$$

For $q^* \in H^*$, we have

$$\chi^*_{B_\delta}(q^*) = \sup_{q \in H} \left( \langle q, q^*\rangle_{H, H^*} - \chi_{B_\delta}(q) \right),$$

that is

$$\chi^*_{B_\delta}(q^*) = \sup_{q \in H,\, \|q - g^\delta\|_H \leq \delta} \langle q, q^*\rangle_{H, H^*} = \langle g^\delta, q^*\rangle_{H, H^*} + \delta \sup_{q \in H,\, \|q\|_H \leq 1} \langle q, q^*\rangle_{H, H^*}$$

$$= \langle g^\delta, q^*\rangle_{H^*, H} + \delta \|q^*\|_{H^*}.$$

The problem $(\mathscr{P}^*)$ is hence

$$(\mathscr{P}^*) \quad \sup_{q^* \in H^*} -\Phi^*(0, q^*) = \sup_{q^* \in H} \left( -\frac{1}{2}\|A^*q^*\|_{V^*}^2 - \delta \|q^*\|_{H^*} + \langle g^\delta, q^*\rangle_{H, H^*} \right)$$

and we finally obtain problem (44). $\square$

By applying Theorem 5.1 using the identification $V^* = V$ and $H^* = H$ by the Riesz theorem, we can now emphasize the link between the solution to the primal problem ($\mathscr{P}$) and the solutions to the dual problem ($\mathscr{P}^*$).

**Proposition 5.3.** *If the assumption (1) is satisfied, the problem ($\mathscr{P}^*$) has solutions. In addition, the solution $u^\delta$ to the primal problem ($\mathscr{P}$) given by (42) and the solutions $p^\delta$ to the dual problem ($\mathscr{P}^*$) given by (44) are related to each other by $u^\delta = A^* p^\delta$.*

*Proof.* Let us check that we satisfy the assumptions of theorem 5.1. Clearly $\Phi$ given by (43) is a convex function of $(u, q)$ and $\inf(\mathscr{P}) < +\infty$ since problem ($\mathscr{P}$) has a (unique) solution. Let us use the decomposition $g^\delta = g_\parallel^\delta + g_\perp^\delta$, with $g_\parallel^\delta \in \overline{\text{Range } A}$ and $g_\perp^\delta \in (\text{Range } A)^\perp$. In virtue of assumption (1), we have $\delta_\perp = \|g_\perp^\delta\|_H < \delta$. Let us choose $u_0 \in V$ such that $\|Au_0 - g_\parallel^\delta\|_H \leq (\delta - \delta_\perp)/2$. For any $q \in H$ such that $\|q\|_H \leq (\delta - \delta_\perp)/2$, we have

$$\|Au_0 - q - g^\delta\|_H \leq \|Au_0 - g_\parallel^\delta\|_H + \|g_\perp^\delta\|_H + \|q\|_H \leq \delta,$$

that is $q \to \Phi(u_0, q) = L(u_0) < +\infty$ is constant in a neighborhood of point 0. We can then apply Theorem 5.1. In particular, it implies that ($\mathscr{P}^*$) has solutions and that

$$\inf(\mathscr{P}) = \sup(\mathscr{P}^*) < +\infty.$$

From now on we identify $V^*$ and $H^*$ with $V$ and $H$, respectively. Let $p^\delta$ be a solution of ($\mathscr{P}^*$), the above relationship implies that

$$\frac{1}{2}\|u^\delta\|_V^2 = -\frac{1}{2}\|A^* p^\delta\|_V^2 - \delta \|p^\delta\|_H + (g^\delta, p^\delta)_H. \tag{45}$$

We already know from the proof of Theorem 2.2 in the case when $P = 0$ that

$$\|A^* p^\delta\|_V^2 + \delta \|p^\delta\|_H - (g^\delta, p^\delta)_H = 0,$$

which together with (45) implies that

$$\frac{1}{2}\|u^\delta\|_V^2 = \frac{1}{2}\|A^* p^\delta\|_V^2.$$

Since we also have

$$\|A(A^* p^\delta) - g^\delta\|_H = \delta,$$

we conclude that $A^* p^\delta \in V$ solves the primal problem ($\mathscr{P}$), that is $A^* p^\delta = u^\delta$. □

**Remark 5.2.** *It is remarkable that the primal problem ($\mathscr{P}$) is more difficult to solve in practice than the dual problem ($\mathscr{P}^*$), which is unconstrained.*

**Remark 5.3.** *The solution $u^\delta$ to the primal problem ($\mathscr{P}$) coincides with the Tikhonov solution $u_\varepsilon^\delta$ to problem (2) associated with the Morozov value $\varepsilon(\delta)$.*